\documentclass[10pt]{amsart}
\usepackage{amsmath,amsfonts,amssymb,amsthm}
\usepackage[alphabetic]{amsrefs}
\usepackage{url}
\usepackage{graphicx,color}

\newtheorem{theorem}{Theorem}[section]
\newtheorem{lemma}[theorem]{Lemma}

\newtheorem{prop}[theorem]{Proposition}

\newtheorem*{theorem*}{Theorem}
\newtheorem*{corollary}{Corollary}
\theoremstyle{definition}
\newtheorem{definition}{Definition}[section]

\theoremstyle{remark}

\usepackage{caption}
\begin{document}

\title[\ ]{}

\begin{center}
\uppercase{\textbf{Minimal disc diagrams of $5/9$-simplicial complexes}}
\footnote{This work was partially supported by the grant $346300$ for IMPAN from
the Simons Foundation and the matching $2015-2019$ Polish MNiSW fund.}\\
\vspace{0.5cm}
\end{center}

\author[\ ]{
Ioana-Claudia Laz\u{a}r\\
Politehnica University of Timi\c{s}oara, Dept. of Mathematics,\\
Victoriei Square $2$, $300006$-Timi\c{s}oara, Romania\\
E-mail address: ioana.lazar@upt.ro}

\date{}

\begin{abstract}

We introduce and study a local combinatorial condition, called the $5/9$-condition, on a simplicial complex, implying Gromov hyperbolicity of its universal cover.
We hereby give an application of another combinatorial condition, called $8$-location, introduced by Damian Osajda. Along the way we prove the minimal filling diagrams lemma for $5/9$-complexes.

\hspace{0 mm} \textbf{2010 Mathematics Subject Classification}:
05C99, 05C75.

\hspace{0 mm} \textbf{Keywords}: $5/9$-condition, $8$-location, minimal disc diagram, nondege\-nerate map.
\end{abstract}

\pagestyle{myheadings}

\markboth{}{}

  \vspace{-10pt}

\maketitle

\section{Introduction}

Curvature can be expressed both in metric and combinatorial terms. On the metric side, one can refer to nonpositively curved in the sense of Aleksandrov and Gromov,
i.e. by comparing small triangles in the space with triangles in the Euclidean plane. Such triangles must satisfy the CAT(0) inequality. On the combinatorial side, one can express curvature
using a condition, called local $6$-largeness which was introduced independently by Chepoi \cite{Ch} (under the name of bridged complexes), Januszkiewicz-{\' S}wi{\c a}tkowski \cite{JS1} and Haglund \cite{Hag}. A flag simplicial complex is
\emph{locally $k$-large} if its links do not contain essential loops of length less than $k$.
In particular, simply connected locally $7$-large simplicial complexes, i.e.\ \emph{$7$-systolic} complexes, are Gromov hyperbolic (see \cite{JS1}).

The theory of \emph{$7$-systolic groups}, that is groups acting geometrically on $7$-systolic complexes, allowed to provide examples of highly dimensional Gromov hyperbolic groups (see \cites{JS0,JS1,O-chcg,OS}). Some very restrictive limitations
are known, however, for groups acting geometrically on CAT(-1) cubical complexes or on $7$-systolic complexes. For example, $7$-systolic groups are in a sense `asymptotically hereditarily aspherical', i.e.\ asymptotically they can not contain
essential spheres. This yields in particular that such groups are not fundamental groups of negatively curved manifolds of dimension
above two; see e.g.\ \cites{JS2,O-ci,O-ib,OS, Gom,O-ns}.

In \cites{O-sdn, ChOs,BCCGO,ChaCHO} other conditions
of this type are studied. They form a way of unifying CAT(0) cubical and systolic theories.
Osajda introduced in \cite{O-8loc} another local combinatorial condition called \emph{$8$-location}. Under the additional hypothesis
of local $5$-largeness, this condition turned out to be useful in finding a new solution to Thurston's
problem about hyperbolicity of some $3$-manifolds.

In \cite{L-8loc} we study of a version of $8$-location, suggested in \cite[Subsection 5.1]{O-8loc}.
This $8$-location says that homotopically trivial loops of length at most $8$ admit filling diagrams with one internal vertex.
However, in the new $8$-location essential $4$-loops are allowed.
In \cite{L-8loc} (Theorem $4.3$) it is shown that this local combinatorial condition is a negative-curvature-type condition. Namely, we prove that simply connected,
$8$-located simplicial complexes are Gromov hyperbolic. In the current paper we give an application to this result.
We introduce another combinatorial curvature condition, called the $5/9$-condition, and we show that the complexes which fulfill it, are also Gromov hyperbolic. In particular, $5/9$-complexes satisfy a linear isoperimetric inequality (see \cite{BH}, chapter $III.H$, page $417$ and page $419$; \cite{L-8loc2}).
The proof relies on the minimal filling diagrams lemma for $5/9$-complexes.
Namely, we prove that the minimal disc in the filling diagram associated to a cycle in a simply connected $5/9$-complex, is itself a $5/9$-complex (Proposition \ref{3.3}).
Further we show that a $5/9$-disc, is $8$-located (Theorem \ref{5.1}). A simply connected $5/9$-complex is therefore itself $8$-located and hence,
by \cite{L-8loc}, Gromov hyperbolic (Corollary $5.1$).
Our proof for the minimal filling diagrams lemma is similar to the one given in \cite{JS1} (Lemma $1.6$ and Lemma $1.7$) and in \cite{Pr} (Theorem $2.7$) for systolic complexes.

The paper's main result is already mentioned without proof in \cite{O-8loc} (Section $5.3$) as an application concerning
some weakly systolic complexes and groups.
Because the $5/9$-condition implies the hyperbolicity of weakly systolic complexes, it
is more general than the well studied
conditions for hyperbolicity: local $7$-largeness for $7$-systolic complexes (see \cite{JS1}, Section $2$), and SD$_{2}^{\ast}$(7) for weakly systolic complexes (see \cite{O-sdn}, Section $7$).

\medskip

\textbf{Acknowledgements.} The author is grateful to Damian Osajda for introducing her to the subject, posing the problem and helpful remarks. The author would also like to thank Djordje Baralic for hospitality and useful discussions during her visit at the Mathematical Institute of the Serbian Academy of Sciences and Arts, Belgrade in December $2017$.

\section{Preliminaries}

Let $X$ be a simplicial complex. We denote by $X^{(k)}$ the $k$-skeleton of
$X, 0 \leq k < \dim X$. A subcomplex $L$ in $X$ is called \emph{full} (in $X$) if any simplex of $X$ spanned by a set of vertices in $L$, is a simplex of $L$. For a set
$A = \{ v_{1}, ..., v_{k} \}$ of vertices of $X$, by $\langle  A \rangle$ or by $\langle  v_{1}, ..., v_{k} \rangle$ we denote the \emph{span} of $A$, i.e. the
smallest full  subcomplex of $X$ that
contains $A$. We write $v \sim v'$ if $\langle  v,v' \rangle \in X$. We write $v \nsim v'$ if $\langle  v,v' \rangle \notin X$.
 We call $X$ {\it flag} if any finite set of vertices which are pairwise connected by
edges of $X$, spans a simplex of $X$.

A {\it cycle} ({\it loop}) $\gamma$ in $X$ is a subcomplex of $X$ isomorphic to a triangulation of $S^{1}$.
A $k$-\emph{wheel} in $X$ $(v_{0}; v_{1}, ..., v_{k})$ (where $v_{i}, 0 \leq i \leq k$
are vertices of $X$) is a subcomplex of $X$ such that $(v_{1}, ..., v_{k})$ is a full cycle and $v_{0} \sim v_{1}, ..., v_{k}$.
The \emph{length} of $\gamma$ (denoted by $|\gamma|$) is the number of edges of $\gamma$.
If $v_{i} \sim v_{i+j}, 1 \leq i \leq k, 2 \leq j \leq k-i$, then we call $\langle v_{i}, v_{i+j} \rangle$ a \emph{diagonal} of the cycle $\gamma$.

We define the \emph{metric} on the $0$-skeleton of $X$ as the number of edges in the shortest $1$-skeleton path joining two given vertices and we denote it by $d$.
We call two vertices $v,w$ of $X$ \emph{neighbors} if $d(v,w) = 1$. A \emph{ball (sphere)}
$B_{i}(v,X)$ ($S_{i}(v,X)$) of radius $i$ around some vertex $v$ is a full subcomplex of $X$ spanned by vertices at distance at most $i$ (at distance $i$) from $v$.

\begin{definition}\label{2.1}
A simplicial complex is $m$-\emph{located} if it is flag and every full homotopically trivial loop of length at most $m$ is contained in a $1$-ball.
\end{definition}

The \emph{link} of $X$ at $\sigma$, denoted $X_{\sigma}$, is the subcomplex of $X$ consisting of all simplices of $X$ which are disjoint from $\sigma$ and which, together
with $\sigma$, span a simplex of $X$.  A \emph{full cycle} in $X$ is a cycle that is full as subcomplex of $X$. We call a flag simplicial complex $k$\emph{-large} if there are no
full $j$-cycles in $X$, for $j<k$. We say $X$ is \emph{locally} $k$\emph{-large} if all its links are $k$-large.
We call a vertex of $X$ \emph{$k$-large} if its link is $k$-large.

\begin{definition}\label{2.2}
We say that a flag simplicial complex satisfies the \emph{$5/9$-condition}, or that it is a \emph{$5/9$-complex}, if it satisfies the following three conditions:
\begin{description}
\item[($5/9$)]
every vertex adjacent to a non-$5$-large (but $4$-large) vertex is $9$-large;
\item[($6/8$)]
every vertex adjacent to a non-$6$-large (but $5$-large) vertex is $8$-large;
\item[($7/7$)]
every vertex adjacent to a non-$7$-large (but $6$-large) vertex is $7$-large.
\end{description}
\end{definition}

Further we give an example of a complex which does not satisfy the $5/9$-condition. Let $X$ be a flag simplicial $2$-complex obtained by taking a square grid $S$ and adding, for each square of $S$, a vertex adjacent to the four vertices of such square. Let $\sigma_{1} = \langle v_{1},v_{2},v_{4},v_{3} \rangle$, $\sigma_{2} = \langle v_{2},v_{7},v_{8},v_{4} \rangle$, $\sigma_{3} = \langle v_{3},v_{4},v_{6},v_{5} \rangle$ and $\sigma_{4} = \langle v_{4},v_{8},v_{9},v_{6} \rangle$ be four squares of $S$. Let $w_{i} \in X$ be vertices of $X$ such that $w_{i}$ is adjacent to the vertices of $\sigma_{i}$, $1 \leq i \leq 4$. Note that the vertices $w_{1}$ and $v_{4}$ are neighbors and that $w_{1}$ is $4$-large (but not $5$-large) while $v_{4}$ is $8$-large. Because $v_{4}$ is not $9$-large, the condition ($5/9$) in the definition of a $5/9$-complex is not fulfilled. So $X$ is not a $5/9$-complex.

\begin{figure}[h]
    \begin{center}
       \includegraphics[height=3cm]{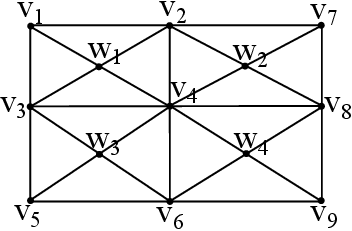}
      \caption{}
    \end{center}
\end{figure}

\begin{definition}\label{2.3}
A \emph{simplicial map} $f : X \rightarrow Y$ between simplicial complexes $X$ and $Y$ is a map
which sends vertices to vertices, and whenever vertices $v_{0}, ..., v_{k} \in X$ span a simplex
$\sigma$ of $X$ then their images span a simplex $\tau$ of $Y$ and we have $f(\sigma) = \tau$.
Therefore a simplicial map is determined by its values on the vertex set of $X$. A
simplicial map is called \emph{nondegenerate} if it is injective on each simplex.
\end{definition}

\begin{definition}\label{2.4}
Let $\gamma$ be a cycle in $X$. A \emph{filling diagram} for $\gamma$ is a simplicial map $f : D \rightarrow X$ where $D$ is a triangulated $2$-disc, and
$f | _{\partial D}$ maps $\partial D$ isomorphically onto $\gamma$. We denote a filling diagram for $\gamma$ by $(D,f)$ and we say it is:

$\bullet$ \emph{minimal} if $D$ consists of the least possible number of $2$-simplices among filling diagrams for $\gamma$;

$\bullet$ \emph{nondegenerate} if $f$ is a nondegenerate map;

$\bullet$ \emph{locally} $k$-\emph{large} if $D$ is a locally $k$-large simplicial complex.

\end{definition}

The disc $D$ in the definition above is locally $k$-large only if it is simplicial.
Also $D$ is contractible; in particular it is simply connected.

\begin{lemma}\label{2.5}
Let $X$ be a simplicial complex and let $\gamma$ be a homotopically trivial loop in $X$. Then:
\begin{enumerate}
\item there exists a filling diagram $(D, f)$ for $\gamma$ (see \cite{Ch} - Lemma $5.1$, \cite{JS1} - Lemma $1.6$ and \cite{Pr} - Theorem $2.7$);
\item any minimal filling diagram for $\gamma$ is simplicial and nondegenerate (see \cite{Ch} - Lemma $5.1$, \cite{JS1} - Lemma $1.6$, Lemma $1.7$ and \cite{Pr} - Theorem $2.7$).

\end{enumerate}
\end{lemma}

\begin{lemma}\label{2.6}
Let $X$ be a simplicial complex and let $\gamma$ be a homotopically trivial loop in $X$. Let $(D,f)$ be a minimal filling diagram for $\gamma$. Then adjacent $2$-simplices of $D$ have distinct images under $f$ (see \cite{Ch} - Lemma $5.1$).

\end{lemma}

\section{Minimal filling diagrams lemma for $5/9$-complexes}

In this section we prove the minimal filling diagrams lemma for $5/9$--complexes.
We start with a useful lemma.

\begin{lemma}\label{3.1}
Let $X$ be a simplicial complex and let $\gamma$ be a homotopically trivial loop in $X$. Let $(D, f)$ be a minimal filling diagram for $\gamma$.
Let $v$ be an interior vertex of $D$ which is $k$-large, $4 \leq k \leq 8$. Then the vertex $f(v)$ differs from any vertex in $X_{f(v)}$. Moreover, the map $f$ is injective on $X_{f(v)}$.
\end{lemma}

\begin{proof}

Let $D_{v} = (v_{1}, v_{2}, ..., v_{k}), 4 \leq k \leq 8$.
Because $(D,f)$ is a minimal filling diagram, the map $f$ is simplicial and nondegenerate. Therefore, since in $D$ there are simplices $v, v_{j}$, $\langle v, v_{j} \rangle$, $1 \leq j \leq k$, $\langle v_{j-1}, v_{j} \rangle$, $2 \leq j \leq k$, in $X$ there are simplices $f(v), f(v_{j})$, $\langle f(v), f(v_{j}) \rangle$, $1 \leq j \leq k$, $\langle f(v_{j-1}), f(v_{j}) \rangle$, $2 \leq j \leq k$.
Lemma \ref{2.6} implies that adjacent $2$-simplices of $D$ have distinct images under $f$.
Hence $f(v) \neq f(v_{i\, mod\, k\; +\; 1}), 1 \leq i \leq k$ and $f(v_{i\, mod\, k\; +\; 1}) \neq f(v_{(i+2)\, mod\, k\; +\; 1}), 1 \leq i \leq k$.

We show further that $f(v_{j\, mod\, k\; +\; 1}) \neq f(v_{(j+3)\, mod\, k\; +\; 1})$, $1 \leq j \leq k$.
Suppose by contradiction there exists $i$ such that $f(v_{i\, mod\, k\; +\; 1}) = f(v_{(i+3)\, mod\, k\; +\; 1})$, $1 \leq i \leq k$.
We choose a filling diagram $(D',f')$ for $\gamma$ such that in $D'$ we have
$v_{i\, mod\, k\; +\; 1} \sim v_{(i+j)\, mod\, k\; +\; 1}$, $2 \leq j \leq 4$.
We triangulate $D'$ with the same simplices like $D$ except for the triangles
$\langle v, v_{(i+j)\, mod\, k\; +\; 1}, v_{(i+j+1)\, mod\, k\; +\; 1} \rangle, 0 \leq j \leq 3$ in $D$ which are replaced in $D'$ by the triangles
$\langle v, v_{i\, mod\, k\; +\; 1}, v_{(i+4)\, mod\, k\; +\; 1} \rangle,$ $\langle v_{i\, mod\, k\; +\; 1}, v_{(i+j)\, mod\, k\; +\; 1}, v_{(i+j+1)\, mod\, k\; +\; 1} \rangle,$ $1 \leq j \leq 3$.
We define $f'$ such that it coincides with $f$ on all simplices which are common to $D$ and $D'$.
We define $f'$ such that $f'(v_{i\, mod\, k\; +\; 1}) = f'(v_{(i+3)\, mod\, k\; +\; 1}) = f(v_{i\, mod\, k\; +\; 1})$,\\
$f'(\langle v_{i\, mod\, k\; +\; 1}, v_{(i+3)\, mod\, k\; +\; 1} \rangle) = f(v_{i\, mod\, k\; +\; 1})$.
As argued above $f(v_{j\, mod\, k\; +\; 1}) \neq f(v_{(j+2)\, mod\, k\; +\; 1}), 1 \leq j \leq k$.
Since in $X$ we have $f(v_{(i+3)\, mod\, k\; +\; 1})$ $\sim$ \\ $f(v_{(i+2)\, mod\, k\; +\; 1})$, we may define $f'$ such that $f'(v_{(i+3)\, mod\, k\; +\; 1})$ $\sim$ $f'(v_{(i+2)\, mod\, k\; +\; 1})$.
Then because $f'(v_{i\, mod\, k\; +\; 1})$ $=$ $f'(v_{(i+3)\, mod\, k\; +\; 1})$, we have $f'(v_{i\, mod\, k\; +\; 1})$ $\sim$ \\ $f'(v_{(i+2)\, mod\, k\; +\; 1})$.
We define $f'$ such that $f'(\langle v_{i\, mod\, k\; +\; 1}, v_{(i+2)\, mod\, k\; +\; 1} \rangle)$ $=$ \\ $\langle f'(v_{i\, mod\, k\; +\; 1}), f'(v_{(i+2)\, mod\, k\; +\; 1}) \rangle$ $=$ $\langle f(v_{i\, mod\, k\; +\; 1}), f(v_{(i+2)\, mod\, k\; +\; 1}) \rangle$.
One can similarly show that we can define $f'$ such that $f'(\langle v_{i\, mod\, k\; +\; 1}, v_{(i+4)\, mod\, k\; +\; 1} \rangle) = \langle f'(v_{i\, mod\, k\; +\; 1}), f'(v_{(i+4)\, mod\, k\; +\; 1}) \rangle = \langle f(v_{i\, mod\, k\; +\; 1}), f(v_{(i+4)\, mod\, k\; +\; 1}) \rangle$.
We define $f'$ such that $f'(\langle v, v_{i\, mod\, k\; +\; 1}, v_{(i+4)\, mod\, k\; +\; 1} \rangle)$ $=$ \\ $\langle f(v), f(v_{i\, mod\, k\; +\; 1}),$
$f(v_{(i+4)\, mod\, k\; +\; 1})$ $\rangle$,\\
$f'(\langle v_{i\, mod\, k\; +\; 1}, v_{(i+j)\, mod\, k\; +\; 1}, v_{(i+j+1)\, mod\, k\; +\; 1} \rangle)$ $=$ \\ $\langle f(v_{i\, mod\, k\; +\; 1}), f(v_{(i+j)\, mod\, k\; +\; 1}), f(v_{(i+j+1)\, mod\, k\; +\; 1}) \rangle, 1 \leq j \leq 3$.
Hence, since $f$ is simplicial, $f'$ is also simplicial.
So $(D', f')$ is indeed a filling diagram for $\gamma$.
Note that $D$ and $D'$ have the same area.
Therefore $D'$ has minimal area. Then Lemma \ref{2.1} implies that the map $f'$ is nondegenerate.
But since $f'(v_{i\, mod\, k\; +\; 1}) = f'(v_{(i+3)\, mod\, k\; +\; 1})$, $f'$ is degenerate.
Because we have reached a contradiction, $f(v_{i\, mod\, k\; +\; 1}) \neq f(v_{(i+3)\, mod\, k\; +\; 1}),$ $1 \leq i \leq k$.

To finish the proof in case $k = 8$, we still need to show that the vertex $f(v_{1})$ differs from the vertex $f(v_{5})$.
Suppose by contradiction $f(v_{1}) = f(v_{5})$. We choose a filling diagram $(D',f')$ for $\gamma$ such that $D'$ is triangulated with the same simplices like $D$.
As argued above $f(v) \neq f(v_{i\, mod\, k\; +\; 1})$, $f(v_{i\, mod\, k\; +\; 1}) \neq f(v_{(i+2)\, mod\, k\; +\; 1})$ and $f(v_{i\, mod\, k\; +\; 1}) \neq f(v_{(i+3)\, mod\, k\; +\; 1}), 1 \leq i \leq 8$.
We define $f'$ such that it coincides with $f$ on all simplices which are common to $D$ and $D'$.
We define $f'$ such that $f'(v_{1}) = f'(v_{5}) = f(v_{1})$, $f'(\langle v, v_{5} \rangle) = \langle f(v), f(v_{1}) \rangle$, $f'(\langle v_{i},v_{1} \rangle) = \langle f(v_{i}),f(v_{1}) \rangle$, $f'(\langle v,v_{i},v_{1} \rangle) = \langle f(v),f(v_{i}),f(v_{1}) \rangle$, $i \in \{2,8\}$, $f'(\langle v_{i},v_{5} \rangle) = \langle f(v_{i}),f(v_{1}) \rangle$, \\ $f'(\langle v,v_{i},v_{5} \rangle) = \langle f(v),f(v_{i}),f(v_{1}) \rangle$, $i \in \{4,6\}$.
Because $f$ is simplicial, $f'$ is also simplicial. So $(D', f')$ is indeed a filling diagram for $\gamma$.
Note that the discs $D$ and $D'$ have the same area.
Hence $D'$ has minimal area and therefore $f'$ is nondegenerate.
But, since $f'(v_{1}) = f'(v_{5})$, the map $f'$ is degenerate.
Because we have reached a contradiction, $f(v_{1}) \neq f(v_{5})$.

In conclusion $f(v)$ differs from any vertex in $X_{f(v)}$ and $f$ is injective on $X_{f(v)}$.

\end{proof}

The following lemma will be used frequently.

\begin{lemma}\label{3.2}
Let $X$ be a $5/9$-complex and let $(k,k')$ be one of the pairs $(4,9)$, $(5,8)$, $(6,7)$. Let $v$ and $w$ be adjacent vertices of $X$. The vertex $v$ is $k$-large but not $(k+1)$-large and the vertex $w$ is $k''$-large, $k'' < k'$. Let $\gamma$ be a cycle in $X_{w}$. Then $\gamma$ is not full in $X$.
\end{lemma}

\begin{proof}
Let $(k,k') = (4,9)$. Because $X$ is a $5/9$-complex, any neighbor of $v$ must be $9$-large. Therefore, since the length of the cycle $\gamma$ is at most $8$, $\gamma$ is not full in $X$. We argue similarly for the cases $(k,k') \in \{(5,8), (6,7)\}$.
\end{proof}

We prove further the minimal filling diagrams lemma for $5/9$-complexes.

\begin{prop}[$5/9$-disc diagrams]\label{3.3}
Let $X$ be a $5/9$-complex and let $\gamma$ be a homotopically trivial loop in $X$. Then for any minimal filling diagram $(D, f)$ for $\gamma$, $D$ is a $5/9$-complex.

\end{prop}

\begin{proof}

We will prove that $D$ is a $5/9$-complex as follows.
Condition ($5/9$) is shown in the cases $1 - 5$ below. Condition ($6/8$) is shown in the cases $1 - 4$ below. Condition ($7/7$) is shown in the cases $1 - 3$ below.

Let $(k,k') \in \{(4,9)$, $(5,8)$, $(6,7) \}$.
Let $v$ be an interior vertex of $D$ which is $k$-large but not $(k+1)$-large. So $D_{v}$ is a full cycle of length $k$ of $D$. We will show that any neighbor of $v$ is $k'$-large. Let $w$ be an interior vertex of $D$ adjacent to $v$. Then $w \in D_{v}$ and $v \in D_{w}$.

Because $D$ has minimal area,
Lemma \ref{2.1} implies that $f$ is simplicial and nondegenerate. Hence $f(v)$ and $f(w)$ are adjacent. Let $a \in \{v,w\}$.
According to Lemma \ref{3.1}, the map $f$ is injective on $X_{f(a)}$.
So no two vertices in $X_{f(a)}$ coincide with each other.
Lemma \ref{3.1} also implies that the vertex $f(a)$ differs from any vertex in $X_{f(a)}$.

Our subsequent proof is to show that $w$ is $k''$-large for all values of $k'', k'' \leq k'$. Because $D$ is flat, the flagness of $D$ implies that $k'' \geq 4$. We distinguish five cases depending on the value of $k''$. Namely $4 \leq k'' \leq 8$.

\begin{enumerate}

\item Case $1$: $k'' = 4$.

Let $D_{v} = (v_{1}, w, v_{3}, ..., v_{k})$, $4 \leq k \leq 6$ and let $D_{w} = (v_{1}, v_{2}, v_{3}, v)$. Because $f$ is simplicial and nondegenerate, there is a $k$-cycle $(f(v_{1}), f(w), f(v_{3}),$ $..., f(v_{k}))$, $4 \leq k \leq 6$ in $X_{f(v)}$ and there is a $4$-cycle $\alpha = (f(v_{1}), f(v_{2}), f(v_{3}),$ $f(v))$ in
$X_{f(w)}$. Lemma \ref{3.2} implies that $\alpha$ is not full in $X$. We distinguish the cases $f(v_{1}) \sim f(v_{3})$, $f(v) \sim f(v_{2})$ which are analyzed below.

\begin{figure}[h]
    \begin{center}
       \includegraphics[height=4cm]{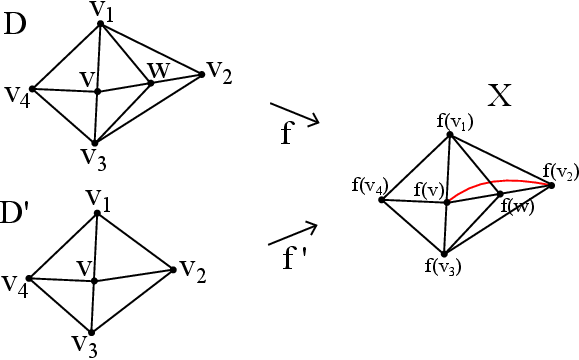}
      \caption{Case $1$, $k = 4$, $f(v) \sim f(v_{2})$}
    \end{center}
\end{figure}

\begin{figure}[h]
    \begin{center}
       \includegraphics[height=4cm]{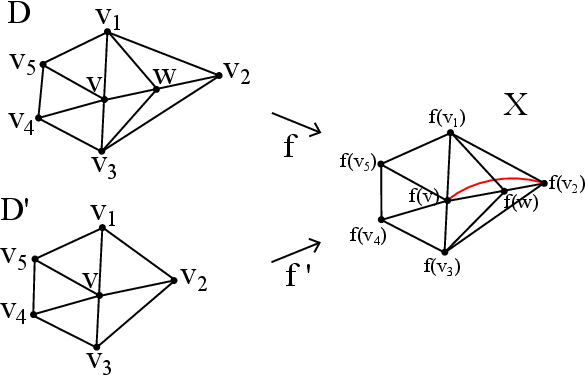}
      \caption{Case $1$, $k = 5$, $f(v) \sim f(v_{2})$}
    \end{center}
\end{figure}

\begin{figure}[h]
    \begin{center}
       \includegraphics[height=4cm]{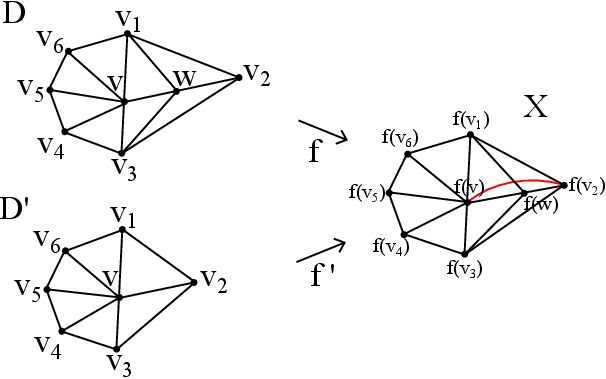}
      \caption{Case $1$, $k = 6$, $f(v) \sim f(v_{2})$}
    \end{center}
\end{figure}

a. Suppose $f(v_{1}) \sim f(v_{3})$.

We consider a filling diagram $(D',f')$ for $\gamma$ such that $D'$ is triangulated with the same simplices like $D$ except for the triangles $\langle v,w,v_{1} \rangle$, $\langle v,w,v_{3} \rangle$, $\langle v,v_{i},v_{i+1}\rangle$, $1 \leq i \leq 2$ in $D$ which are replaced in $D'$ by the triangles $\langle v,v_{1},v_{3} \rangle$, $\langle v_{1},v_{2},v_{3} \rangle$. We define $f'$ such that it coincides with $f$ on all simplices which are common to $D$ and $D'$. We define $f'$ such that $f'(\langle v_{1},v_{3} \rangle) = \langle f(v_{1}),f(v_{3}) \rangle$, $f'(\langle v,v_{1},v_{3} \rangle) = \langle f(v),f(v_{1}),f(v_{3}) \rangle$, $f'(\langle v_{1},v_{2},$ $v_{3} \rangle) = \langle f(v_{1}),f(v_{2}),f(v_{3}) \rangle$. Because $f$ is simplicial, $f'$ is also simplicial. So $(D',f')$ is a filling diagram for $\gamma$. Note that the area of $D'$ is less than the area of $D$. Because $D$ has minimal area, we have reached a contradiction.

b. Suppose $f(v) \sim f(v_{2})$. We choose a filling diagram $(D', f')$ for $\gamma$ such that $D'$ is triangulated with the same simplices like $D$ except for the triangles $\langle v, w, v_{1} \rangle, $ $\langle v, w, v_{3} \rangle,$
$ \langle w, v_{i}, v_{i+1} \rangle,$
$1 \leq i \leq 2$ in $D$ which are replaced in $D'$ by the triangles $\langle v, v_{i}, v_{i+1} \rangle,$ $1 \leq i \leq 2$. We define $f'$ such that it coincides with $f$ on all simplices which are common to $D$ and $D'$.
We define $f'$ such that $f'(\langle v, v_{2} \rangle) = \langle f(v), f(v_{2}) \rangle$, $f'(\langle v,v_{i},v_{i+1} \rangle) = \langle f(v),f(v_{i}),f(v_{i+1}) \rangle, 1 \leq i \leq 2$. Because $f$ is simplicial, $f'$ is also simplicial. Then $(D', f')$ is indeed a filling diagram for $\gamma$. The area of $D'$ is less than the area of $D$. Because $D$ has minimal area, we have reached a contradiction.

So $w$ is $5$-large.

\item Case $2$: $k'' = 5$.

\begin{figure}[h]
    \begin{center}
       \includegraphics[height=4cm]{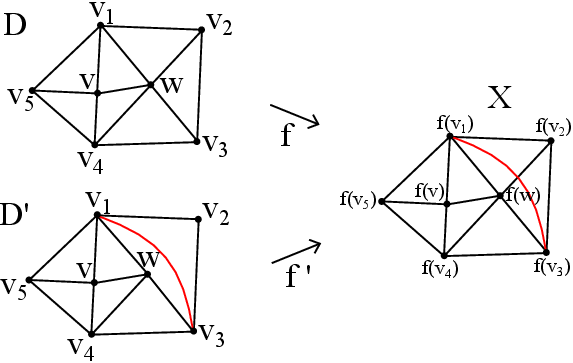}
      \caption{Case $2$, $k = 4$, $f(v_{1}) \sim f(v_{3})$}
    \end{center}
\end{figure}

\begin{figure}[h]
    \begin{center}
       \includegraphics[height=4cm]{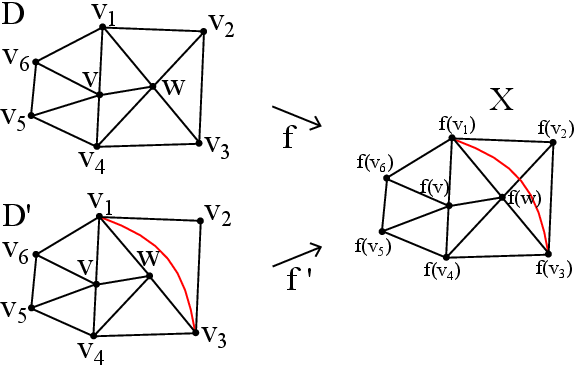}
      \caption{Case $2$, $k = 5$, $f(v_{1}) \sim f(v_{3})$}
    \end{center}
\end{figure}

\begin{figure}[h]
    \begin{center}
       \includegraphics[height=4cm]{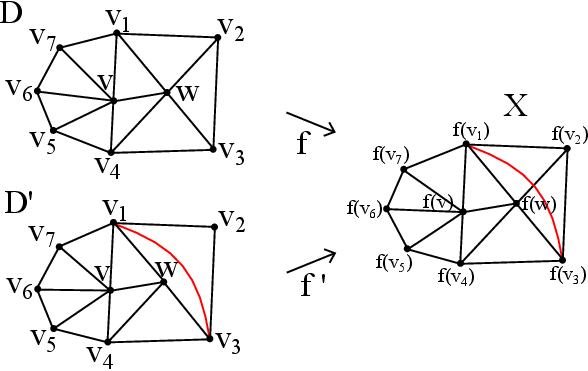}
      \caption{Case $2$, $k = 6$, $f(v_{1}) \sim f(v_{3})$}
    \end{center}
\end{figure}

Let $D_{v} = (v_{1}, w, v_{4}, ..., v_{k+1})$, $4 \leq k \leq 6$ and let $D_{w} = (v_{1}, v_{2}, v_{3}, v_{4}, v)$.
Lemma \ref{3.2} implies that the $5$-cycle
$(f(v_{1}), f(v_{2}), f(v_{3}), f(v_{4}), f(v))$ in $X_{f(w)}$ is not full. We distinguish the following cases: $f(v_{1}) \sim f(v_{i})$, $3 \leq i \leq 4$, $f(v) \sim f(v_{2})$ discussed below.

a. Suppose $f(v_{1}) \sim f(v_{3})$.
We choose a filling diagram $(D', f')$ for $\gamma$ such that $D'$ is triangulated with the same simplices like $D$ except for the triangles $\langle w, v_{i}, v_{i+1} \rangle,$ $1 \leq i \leq 2$ in $D$ which are replaced in
$D'$ by the triangles $\langle w, v_{1}, v_{3} \rangle,$ $\langle v_{1}, v_{2}, v_{3} \rangle$.
We define $f'$ such that it coincides with $f$ on all simplices
which are common to $D$ and $D'$. We define $f'$ such that $f'(\langle v_{1}, v_{3} \rangle) = \langle f(v_{1}), f(v_{3}) \rangle$, $f'(\langle w,v_{1},v_{3} \rangle) = \langle f(w),f(v_{1}),f(v_{3}) \rangle$, $f'(\langle v_{1},v_{2},v_{3} \rangle) =$ $ \langle f(v_{1}),$ $f(v_{2}),f(v_{3}) \rangle$.
Because $f$ is simplicial, so is $f'$. Then $(D', f')$ is indeed a filling diagram for $\gamma$.
The discs $D$ and $D'$ have the same area. So $D'$ has minimal area.
Note that there is a full $4$-cycle $(v_{1}, v_{3}, v_{4}, v)$ in $D'_{w}$. So $D'_{w}$ is $4$-large.
Also there is a full $k$-cycle $(v_{1}, w, v_{4}, ..., v_{k+1})$, $4 \leq k \leq 6$ in $D'_{v}$. Hence $D'_{v}$ is $k$-large,
$4 \leq k \leq 6$. Since the vertices $v$ and $w$ are adjacent, case $1$ implies a contradiction.

b. Suppose $f(v_{1}) \sim f(v_{4})$.

Note that there is a $k$-cycle $(f(v_{1}),f(v_{4}),...,f(v_{k+1}))$ in $X_{f(v)}$, $4 \leq k \leq 6$ and there is a $4$-cycle $\alpha = (f(v_{1}),f(v_{2}),f(v_{3}),f(v_{4}))$ in $X_{f(w)}$. Lemma $\ref{3.2}$ implies that $\alpha$ is not full. Assume w.l.o.g. $f(v_{1}) \sim f(v_{3})$. We consider a filling diagram $(D',f')$ for $\gamma$ such that $D'$ is triangulated with the same simplices like $D$ except for the triangles $\langle v,w,v_{1} \rangle$, $\langle v,w,v_{4} \rangle$, $\langle w,v_{i},v_{i+1} \rangle$, $1 \leq i \leq 3$ in $D$ which are replaced in $D'$ by the triangles $\langle v,v_{1},v_{4} \rangle, \langle v_{1},v_{i},$ $v_{i+1} \rangle, 2 \leq i \leq 3$. We define $f'$ such that it coincides with $f$ on all simplices which are common to $D$ and $D'$. We define $f'$ such that $f'(\langle v_{1},v_{i} \rangle) = \langle f(v_{1}),f(v_{i}) \rangle, 3 \leq i \leq 4$, $f'(\langle v_{1},v_{i},v_{i+1} \rangle) = \langle f(v_{1}),f(v_{i}),$ $f(v_{i+1}) \rangle$, $2 \leq i \leq 3$, $f'(\langle v,v_{1},v_{4} \rangle) = \langle f(v),f(v_{1}),f(v_{4}) \rangle$. Because $f$ is simplicial, $f'$ is also simplicial. So $(D',f')$ is indeed a filling diagram for $\gamma$. Note that the area of $D'$ is less than the area of $D$. So we have reached a contradiction.

c. Suppose $f(v) \sim f(v_{2})$. We choose a filling diagram $(D', f')$ for $\gamma$ such that $D'$ is triangulated with the same simplices like $D$ except for the triangles $\langle w, v, v_{1} \rangle, \langle w, v_{1}, v_{2} \rangle$  in $D$ which are replaced in
$D'$ by the triangles $\langle w, v, v_{2} \rangle,$ $\langle v, v_{1}, v_{2} \rangle$.
We define $f'$ such that it coincides with $f$ on all simplices
which are common to $D$ and $D'$. We define $f'$ such that $f'(\langle v, v_{2} \rangle) = \langle f(v), f(v_{2}) \rangle$, $f'(\langle v,w,v_{2} \rangle) = \langle f(v),f(w),f(v_{2}) \rangle$, $f'(\langle v,v_{1},v_{2} \rangle) = $ $ \langle f(v),$ $f(v_{1}),f(v_{2}) \rangle$.
Because $f$ is simplicial, so is $f'$. Then $(D', f')$ is indeed a filling diagram for $\gamma$.
The discs $D$ and $D'$ have the same area. So $D'$ has minimal area.
Note that there is a full $(k+1)$-cycle $(v_{1},v_{2},w,v_{4}, ..., v_{k+1})$, $4 \leq k \leq 6$ in $D'_{v}$ and there is a full $4$-cycle $(v,v_{2},v_{3},v_{4})$ in $D'_{w}$. Hence $D'_{v}$ is $(k+1)$-large and $D'_{w}$ is $4$-large. If $k=4$, then case $1$, $k=5$ implies a contradiction. If $k=5$, then case $1$, $k=6$ implies a contradiction.
Let $k = 6$. Because $D'$ is minimal, $f'$ is simplicial and nondegenerate. So there is a $4$-cycle $\eta_{1} = (f'(v),f'(v_{2}),$ $f'(v_{3}),f'(v_{4}))$ in $X_{f'(w)}$ and there is a $7$-cycle $\eta_{2} = (f'(v_{1}),f'(v_{2}),f'(w),f'(v_{4}),$ $f'(v_{5}), f'(v_{6}), f'(v_{7}))$ in $X_{f'(v)}$. Lemma \ref{3.2} implies that $\eta_{2}$ is not full. So at least one of the following holds $f'(v_{1}) \sim f'(v_{i})$, $4 \leq i \leq 6$, $f'(v_{1}) \sim f'(w)$. Below we treat some of these cases. The other cases can be treated similarly.


\begin{figure}[h]
    \begin{center}
       \includegraphics[height=4cm]{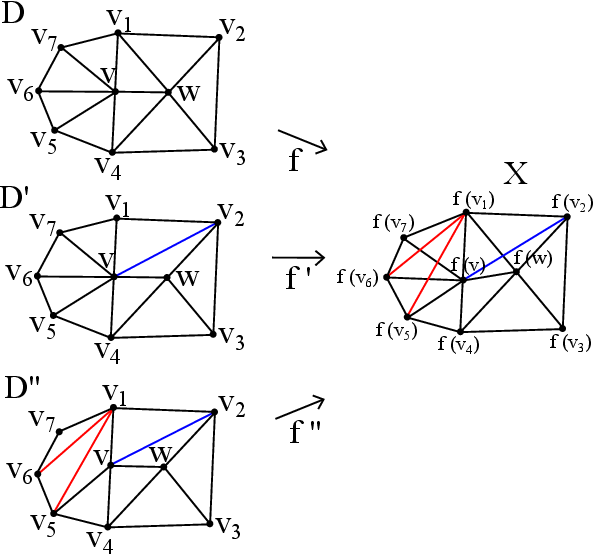}
      \caption{Case $2$, $k=6$, $f(v) \sim f(v_{2}), f'(v_{1}) \sim f'(v_{5})$}
    \end{center}
\end{figure}

$\bullet$ Suppose $f'(v_{1}) \sim f'(v_{4})$. Note that there is a $5$-cycle $\alpha = (f'(v_{1}),$ $f'(v_{4}),$ $f'(v_{5}),f'(v_{6}),f'(v_{7}))$ in $X_{f'(v)}$ and there is a $4$-cycle $(f'(v),f'(v_{2}),$ $f'(v_{3}),$ $f'(v_{4}))$ in $X_{f'(w)}$. By Lemma \ref{3.2}, $\alpha$ is not full. So at least one of the following holds $f'(v_{1}) \sim f'(v_{i})$, $5 \leq i \leq 6$. Assume w.l.o.g. $f'(v_{1}) \sim f'(v_{5})$. The other cases can be treated similarly. By Lemma \ref{3.2}, the $4$-cycle $(f'(v_{1}), f'(v_{5}),$ $f'(v_{6}),f'(v_{7}))$ in $X_{f'(v)}$ is not full. Assume w.l.o.g. $f'(v_{1}) \sim f'(v_{6})$. We choose a filling diagram $(D'',f'')$ for $\gamma$ such that $D''$ is triangulated with the same simplices like $D'$ except for the triangles $\langle v,v_{i},v_{i+1}\rangle, 4 \leq i \leq 6$, $\langle v,v_{1},v_{7}\rangle$ in $D'$ which are replaced  in $D''$ by the triangles $\langle v,v_{1},v_{4}\rangle$, $\langle v_{1},v_{i},v_{i+1}\rangle, 4 \leq i \leq 6$. We define $f''$ such that it coincides with $f'$ on all simplices which are common to $D'$ and $D''$. We define $f''$ such that $f''(\langle v_{1},v_{i}\rangle) = \langle f'(v_{1}), f'(v_{i}) \rangle$, $f''(\langle v_{1},v_{i},v_{i+1} \rangle) = \langle f'(v_{1}), f'(v_{i}), f'(v_{i+1}) \rangle$, $4 \leq i \leq 6$, $f''(\langle v,v_{1},v_{4} \rangle) = \langle f'(v),f'(v_{1}),f'(v_{4}) \rangle$. Because $f'$ is simplicial, $f''$ is also simplicial. So $(D'',f'')$ is indeed a filling diagram for $\gamma$. Note that $D'$ and $D''$ have the same area. Therefore $D''$ has minimal area. There is a full $4$-cycle $(v_{1},v_{2},w,v_{4})$ in $D''_{v}$ and a full $4$-cycle $(v,v_{2},v_{3},v_{4})$ in $D''_{w}$. Hence case $1$, $k = 4$ implies a contradiction.

$\bullet$ Suppose $f'(v_{1}) \sim f'(v_{5})$. By Lemma \ref{3.2}, the $4$-cycle $(f'(v_{1}),f'(v_{5}),$ $f'(v_{6}),f'(v_{7}))$ in $X_{f'(v)}$ is not full. Assume w.l.o.g. $f'(v_{1}) \sim f'(v_{6})$. We choose a filling diagram $(D'',f'')$ for $\gamma$ such that $D''$ is triangulated with the same simplices like $D'$ except for the triangles $\langle v,v_{i},v_{i+1}\rangle, 5 \leq i \leq 6$, $\langle v,v_{1},v_{7}\rangle$ in $D'$ which are replaced  in $D''$ by the triangles $\langle v,v_{1},v_{5}\rangle$, $\langle v_{1},v_{i},v_{i+1}\rangle, 5 \leq i \leq 6$. We define $f''$ such that it coincides with $f'$ on all simplices which are common to $D'$ and $D''$. We define $f''$ such that $f''(\langle v_{1},v_{i}\rangle) = \langle f'(v_{1}), f'(v_{i}) \rangle$, $f''(\langle v_{1},v_{i},v_{i+1} \rangle) = \langle f'(v_{1}), f'(v_{i}), f'(v_{i+1}) \rangle$, $5 \leq i \leq 6$, $f''(\langle v,v_{1},v_{5} \rangle) = \langle f'(v),f'(v_{1}),f'(v_{5}) \rangle$. Because $f'$ is simplicial, $f''$ is also simplicial. So $(D'',f'')$ is indeed a filling diagram for $\gamma$. The area of $D''$ equals the area of $D'$ and therefore $D''$ has minimal area. Note that there is a full $5$-cycle $(v_{1},v_{2},w,v_{4},v_{5})$ in $D''_{v}$, and a full $4$-cycle $(v,v_{2},v_{3},v_{4})$ in $D''_{w}$. So case $1$, $k=5$ implies a contradiction.

$\bullet$ Suppose $f'(v_{1}) \sim f'(v_{6})$. We choose a filling diagram $(D'',f'')$ for $\gamma$ such that $D''$ is triangulated with the same simplices like $D'$ except for the triangles $\langle v,v_{1},v_{7}\rangle$, $\langle v,v_{6},v_{7}\rangle$ in $D'$ are replaced by the triangles $\langle v,v_{1},v_{6}\rangle$, $\langle v_{1},v_{6},v_{7}\rangle$ in $D''$. We define $f''$ such that it coincides with $f'$ on all simplices which are common to $D'$ and $D''$. We define $f''$ such that $f''(\langle v_{1},v_{6}\rangle) = \langle f'(v_{1}), f'(v_{6}) \rangle$, $f''(\langle v,v_{1},v_{6} \rangle) = \langle f'(v),f'(v_{1}),f'(v_{6}) \rangle$, $f''(\langle v_{1},v_{6},v_{7} \rangle) = \langle f'(v_{1}),f'(v_{6}),f'(v_{7}) \rangle$. Because $f'$ is simplicial, $f''$ is also simplicial. So $(D'',f'')$ is indeed a filling diagram for $\gamma$. The area of $D''$ equals the area of $D'$ and therefore $D''$ has minimal area. Note that there is a full $6$-cycle $(v_{1},v_{2},w,v_{4},v_{5},v_{6})$ in $D''_{v}$, and a full $4$-cycle $(v,v_{2},v_{3},v_{4})$ in $D''_{w}$. So case $1$, $k=6$ implies a contradiction.

So $w$ is $6$-large.

\item Case $3$: $k'' = 6$.

Let $D_{v} = (v_{1}, w, v_{5}, ..., v_{k+2})$, $4 \leq k \leq 6$ and let $D_{w} = (v_{1}, v_{2}, v_{3}, v_{4}, v_{5}, v)$.
Lemma \ref{3.2} implies that the $6$-cycle
$(f(v_{1}), f(v_{2}), f(v_{3}), f(v_{4}), f(v_{5}), f(v))$ in $X_{f(w)}$ is not full. We distinguish the cases $f(v_{1}) \sim f(v_{i})$, $3 \leq i \leq 5$, $f(v) \sim f(v_{i})$, $2 \leq i \leq 3$ treated below.

a. Suppose $f(v_{1}) \sim f(v_{3})$. We choose a filling diagram $(D', f')$ for $\gamma$ such that $D'$ is triangulated with the same simplices like $D$ except for the triangles $\langle w, v_{i}, v_{i+1} \rangle, 1 \leq i \leq 2$  in $D$ which are replaced in
$D'$ by the triangles $\langle w, v_{1}, v_{3} \rangle,$ $\langle v_{1}, v_{2}, v_{3} \rangle$.
We define $f'$ such that it coincides with $f$ on all simplices
which are common to $D$ and $D'$. We define $f'$ such that $f'(\langle v_{1}, v_{3} \rangle) = \langle f(v_{1}), f(v_{3}) \rangle$, $f'(\langle w,v_{1},v_{3} \rangle) = \langle f(w),f(v_{1}),f(v_{3}) \rangle$, $f'(\langle v_{1},v_{2},v_{3} \rangle) =$ \\ $\langle f(v_{1}),f(v_{2}),f(v_{3}) \rangle$.
Because $f$ is simplicial, so is $f'$. Then $(D', f')$ is indeed a filling diagram for $\gamma$.
Note that the discs $D$ and $D'$ have the same area. So $D'$ has minimal area. There is a full $k$-cycle $(v_{1},w,v_{5}, ... ,v_{k+2}), 4 \leq k \leq 6$ in $D'_{v}$ and a full $5$-cycle $(v_{1},v_{3},v_{4},v_{5},v)$ in $D'_{w}$. So case $2$ implies a contradiction.

\begin{figure}[h]
    \begin{center}
       \includegraphics[height=4cm]{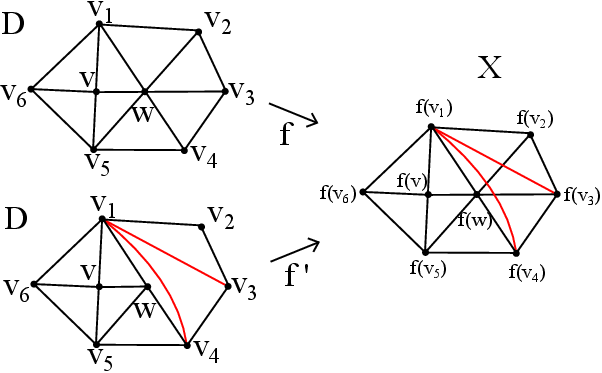}
      \caption{Case $3$, $k = 4$, $f(v_{1}) \sim f(v_{4})$}
    \end{center}
\end{figure}

\begin{figure}[h]
    \begin{center}
       \includegraphics[height=4cm]{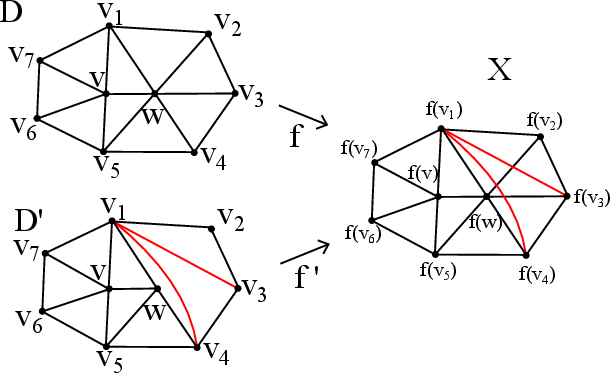}
      \caption{Case $3$, $k = 5$, $f(v_{1}) \sim f(v_{4})$}
    \end{center}
\end{figure}

\begin{figure}[h]
    \begin{center}
       \includegraphics[height=4cm]{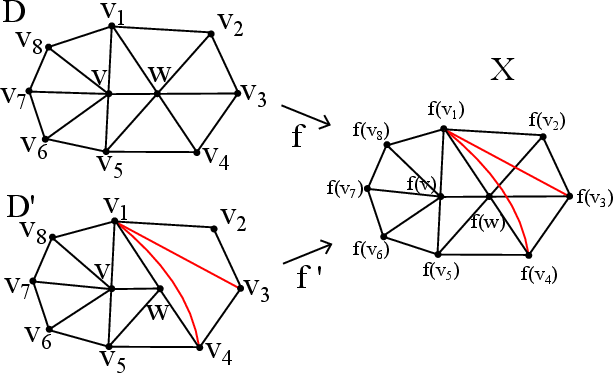}
      \caption{Case $3$, $k = 6$, $f(v_{1}) \sim f(v_{4})$}
    \end{center}
\end{figure}

b. Suppose $f(v_{1}) \sim f(v_{4})$.
Lemma \ref{3.2} implies that the $4$-cycle $(f(v_{1}),$ $ f(v_{2}), f(v_{3}), f(v_{4}))$ in $X_{f(w)}$ is not full. Assume w.l.o.g. $f(v_{1}) \sim f(v_{3})$.
We choose a filling diagram $(D', f')$ for $\gamma$ such that $D'$ is triangulated with the same simplices like $D$
except for the triangles $\langle w, v_{i}, v_{i+1} \rangle, 1 \leq i \leq 3$ in $D$ which are replaced in $D'$ by the triangles
$\langle w, v_{1}, v_{4} \rangle, \langle v_{1}, v_{i}, v_{i+1} \rangle$, $2 \leq i \leq 3$.
 We define $f'$ such that it coincides with $f$ all simplices which are common to $D$ and $D'$. We define $f'$ such that
$f'(\langle v_{1}, v_{i} \rangle) = \langle f(v_{1}), f(v_{i}) \rangle, 3 \leq i \leq 4$, $f'(\langle v_{1},v_{i},v_{i+1} \rangle) = \langle f(v_{1}),f(v_{i}),f(v_{i+1}) \rangle, 2 \leq i \leq 3$, $f'(\langle w,v_{1},v_{4} \rangle) = \langle f(w),f(v_{1}),f(v_{4}) \rangle$. Because $f$ is simplicial, $f'$ is also simplicial. Then $(D', f')$ is indeed a filling diagram for $\gamma$. Note that the discs $D$ and $D'$ have the same area.
So $D'$ has minimal area. There is a full
$4$-cycle $(v_{1}, v_{4}, v_{5}, v)$ in $D'_{w}$ and a full $k$-cycle $(v_{1}, w, v_{5}, ..., v_{k+2})$ in $D'_{v}$, $4 \leq k \leq 6$.
Then case $1$ implies a contradiction.

c. Suppose $f(v_{1}) \sim f(v_{5})$.

Note that there is a $k$-cycle $(f(v_{1}),f(w),f(v_{5}),...,f(v_{k+2}))$ in $X_{f(v)}$, $4 \leq k \leq 6$ and there is a $5$-cycle $\alpha = (f(v_{1}),f(v_{2}),f(v_{3}),f(v_{4}),f(v_{5}))$ in $X_{f(w)}$. Lemma $\ref{3.2}$ implies that $\alpha$ is not full. So at least one of the following holds $f(v_{1}) \sim f(v_{i}), 3 \leq i \leq 4$. Assume w.l.o.g. $f(v_{1}) \sim f(v_{4})$. The other cases can be treated similarly. By Lemma $\ref{3.2}$, the $4$-cycle $(f(v_{1}),f(v_{2}),f(v_{3}),f(v_{4}))$ in $X_{f(w)}$ is not full. Assume w.l.o.g. $f(v_{1}) \sim f(v_{3})$. We consider a filling diagram $(D',f')$ for $\gamma$ such that $D'$ is triangulated with the same simplices like $D$ except for the triangles $\langle v,w,v_{1} \rangle$, $\langle v,w,v_{5} \rangle$, $\langle w,v_{i},v_{i+1} \rangle$, $1 \leq i \leq 4$ in $D$ which are replaced in $D'$ by the triangles $\langle v,v_{1},v_{5} \rangle, \langle v_{1},v_{i},v_{i+1} \rangle, 2 \leq i \leq 4$. We define $f'$ such that it coincides with $f$ on all simplices which are common to $D$ and $D'$. We define $f'$ such that $f'(\langle v_{1},v_{i} \rangle) = \langle f(v_{1}),f(v_{i}) \rangle, 3 \leq i \leq 5$, $f'(\langle v_{1},v_{i},v_{i+1} \rangle) = \langle f(v_{1}),f(v_{i}),$ $f(v_{i+1}) \rangle$, $2 \leq i \leq 4$, $f'(\langle v,v_{1},v_{5} \rangle) = \langle f(v),f(v_{1}),f(v_{5}) \rangle$. Because $f$ is simplicial, $f'$ is also simplicial. So $(D',f')$ is indeed a filling diagram for $\gamma$. Note that the area of $D'$ is less than the area of $D$. Because $D$ has minimal area, we have reached a contradiction.

d. Suppose $f(v) \sim f(v_{2})$. We choose a filling diagram $(D', f')$ for $\gamma$ such that $D'$ is triangulated with the same simplices like $D$ except for the triangles $\langle v,w,v_{1} \rangle$, $\langle w,v_{1},v_{2} \rangle$ in $D$ which are replaced in
$D'$ by the triangles $\langle v,v_{1},v_{2} \rangle, \langle v,w,v_{2} \rangle$.
We define $f'$ such that it coincides with $f$
on all simplices which are common to $D$ and $D'$. We define $f'$ such that $f'(\langle v, v_{2} \rangle) = \langle f(v), f(v_{2}) \rangle$, $f'(\langle v,w,v_{2} \rangle) = \langle f(v),f(w),f(v_{2}) \rangle$, $f'(\langle v,v_{1},v_{2}\rangle) =$ \\ $\langle f(v),f(v_{1}),f(v_{2}) \rangle$.
Because $f$ is simplicial, so is $f'$. Then $(D', f')$ is indeed a filling diagram for $\gamma$.
Note that the discs $D$ and $D'$ have the same area. Then $D'$ has minimal area.
There is a full $(k+1)$-cycle $(v_{1},v_{2},w,v_{5},...,v_{k+2})$, $4 \leq k \leq 6$ in $D'_{v}$, and a full $5$-cycle $(v,v_{2},v_{3},v_{4},v_{5})$ in $D'_{w}$. If $k=4$, then case $2$, $k=5$ implies a contradiction. If $k=5$, we get contradiction with case $2$, $k=6$. Let $k=6$. Because $f'$ is simplicial and nondegenerate, there is $5$-cycle
$\eta_{1} = (f'(v_{2}),f'(v_{3}),f'(v_{4}),f'(v_{5}),f'(v))$ in $X_{f'(w)}$ and there is a $7$-cycle $\eta_{2} = (f'(v_{1}),$ $f'(v_{2}),f'(w),f'(v_{5}),f'(v_{6}),$ $f'(v_{7}),$ $f'(v_{8}))$ in $X_{f'(v)}$. Lemma \ref{3.2} implies that $\eta_{2}$ is not full. So at least one of the following holds $f'(v_{1}) \sim f'(v_{i})$, $5 \leq i \leq 7$, $f'(v_{1}) \sim f'(w)$. Below we treat some of these cases. The other cases can be treated similarly.


\begin{figure}[h]
    \begin{center}
       \includegraphics[height=4cm]{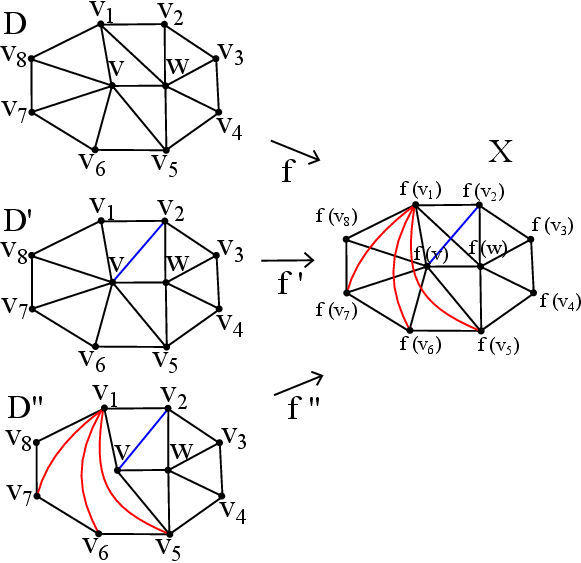}
      \caption{Case $3$, $k=6$, $f(v) \sim f(v_{2}), f'(v_{1}) \sim f'(v_{5})$}
    \end{center}
\end{figure}

$\bullet$ Suppose $f'(v_{1}) \sim f'(v_{5})$. Note that there is a $5$-cycle $\alpha = (f'(v_{1}),$ $f'(v_{5}),$ $f'(v_{6}),$ $f'(v_{7}),f'(v_{8}))$ in $X_{f'(v)}$ and there is a $5$-cycle $(f'(v),f'(v_{2}),$ $f'(v_{3}),$ $f'(v_{4}),f'(v_{5}))$ in $X_{f'(w)}$. By Lemma \ref{3.2}, $\alpha$ is not full. So at least one of the following holds $f'(v_{1}) \sim f'(v_{i})$, $6 \leq i \leq 7$. Assume w.l.o.g. $f'(v_{1}) \sim f'(v_{6})$. The other cases can be treated similarly. By Lemma \ref{3.2}, the $4$-cycle $(f'(v_{1}),f'(v_{6}),$ $f'(v_{7}),f'(v_{8}))$ in $X_{f'(v)}$ is not full. Assume w.l.o.g. $f'(v_{1}) \sim f'(v_{7})$. We consider a filling diagram $(D'',f'')$ for $\gamma$ such that $D''$ is triangulated with the same simplices like $D'$ except for the triangles $\langle v,v_{i},v_{i+1}\rangle, 5 \leq i \leq 7$, $\langle v,v_{1},v_{8}\rangle$ in $D'$ which are replaced in $D''$ by the triangles $\langle v,v_{1},v_{5}\rangle$, $\langle v_{1},v_{i},v_{i+1}\rangle, 5 \leq i \leq 7$. We define $f''$ such that it coincides with $f'$ on all simplices which are common to $D'$ and $D''$. We define $f''$ such that $f''(\langle v_{1},v_{i}\rangle) = \langle f'(v_{1}),f'(v_{i}) \rangle$, $f''(\langle v_{1},v_{i},v_{i+1} \rangle) = \langle f'(v_{1}),f'(v_{i}),f'(v_{i+1}) \rangle$, $5 \leq i \leq 7$, $f''(\langle v,v_{1},v_{5} \rangle) = \langle f'(v),f'(v_{1}),f'(v_{5}) \rangle$. Because $f'$ is simplicial, $f''$ is also simplicial. So $(D'',f'')$ is indeed a filling diagram for $\gamma$. The area of $D''$ equals the area of $D'$ and therefore $D''$ has minimal area. There is a full $4$-cycle $(v_{1},v_{2},w,v_{5})$ in $D''_{v}$, and a full $5$-cycle $(v,v_{2},v_{3},v_{4},v_{5})$ in $D''_{w}$. Then case $2$, $k=4$ (or case $1$, $k=5$) implies a contradiction.

$\bullet$ Suppose $f'(v_{1}) \sim f'(v_{6})$. By Lemma \ref{3.2}, the $4$-cycle $(f'(v_{1}),f'(v_{6}),$ $f'(v_{7}),f'(v_{8}))$ in $X_{f'(v)}$ is not full. Assume w.l.o.g. $f'(v_{1}) \sim f'(v_{7})$. We consider a filling diagram $(D'',f'')$ for $\gamma$ such that $D''$ is triangulated with the same simplices like $D'$ except for the triangles $\langle v,v_{i},v_{i+1}\rangle, 6 \leq i \leq 7$, $\langle v,v_{1},v_{8}\rangle$ in $D'$ which are replaced in $D''$ by the triangles $\langle v,v_{1},v_{6}\rangle$, $\langle v_{1},v_{i},v_{i+1}\rangle, 6 \leq i \leq 7$. We define $f''$ such that it coincides with $f'$ on all simplices which are common to $D'$ and $D''$. We define $f''$ such that $f''(\langle v_{1},v_{i}\rangle) = \langle f'(v_{1}),f'(v_{i}) \rangle$, $f''(\langle v_{1},v_{i},v_{i+1} \rangle) = \langle f'(v_{1}),f'(v_{i}),f'(v_{i+1}) \rangle$, $6 \leq i \leq 7$, $f''(\langle v,v_{1},v_{6}\rangle) = \langle f'(v),f'(v_{1}),f'(v_{6}) \rangle$. Because $f'$ is simplicial, $f''$ is also simplicial. So $(D'',f'')$ is indeed a filling diagram for $\gamma$. The area of $D''$ equals the area of $D'$ and therefore $D''$ has minimal area. There is a full $5$-cycle $(v_{1},v_{2},w,v_{5},v_{6})$ in $D''_{v}$, and a full $5$-cycle $(v,v_{2},v_{3},v_{4},v_{5})$ in $D''_{w}$. By case $2$, $k=5$ we have reached a contradiction.

$\bullet$ Suppose $f'(v_{1}) \sim f'(v_{7})$. We consider a filling diagram $(D'',f'')$ for $\gamma$ such that $D''$ is triangulated with the same simplices like $D'$ except for the triangles $\langle v,v_{1},v_{8}\rangle$, $\langle v,v_{7},v_{8}\rangle$ in $D'$ which are replaced in $D''$ by the triangles $\langle v,v_{1},v_{7}\rangle$, $\langle v_{1},v_{7},v_{8}\rangle$. We define $f''$ such that it coincides with $f'$ on all simplices which are common to $D'$ and $D''$. We define $f''$ such that $f''(\langle v_{1},v_{7}\rangle) = \langle f'(v_{1}),f'(v_{7}) \rangle$, $f''(\langle v,v_{1},v_{7}\rangle) = \langle f'(v),f'(v_{1}),f'(v_{7}) \rangle$, $f''(\langle v_{1},v_{7},v_{8}\rangle) = \langle f'(v_{1}),f'(v_{7}),f'(v_{8}) \rangle$. Because $f'$ is simplicial, $f''$ is also simplicial. So $(D'',f'')$ is indeed a filling diagram for $\gamma$. The area of $D''$ equals the area of $D'$ and therefore $D''$ has minimal area. There is a full $6$-cycle $(v_{1},v_{2},w,v_{5},v_{6},v_{7})$ in $D''_{v}$, and a full $5$-cycle $(v,v_{2},v_{3},v_{4},v_{5})$ in $D''_{w}$. Then case $2$, $k=6$ implies a contradiction.



e. Suppose $f(v) \sim f(v_{3})$. Note that there is a $4$-cycle
$\alpha = (f(v), f(v_{1}),$ $f(v_{2}),$ $f(v_{3}))$ in $X_{f(w)}$ and there is a $k$-cycle
$(f(v_{1}), f(w), f(v_{5}), ... ,f(v_{k+2}))$ in $X_{f(v)}$, $4 \leq k \leq 6$. Lemma \ref{3.2} implies that $\alpha$ is not full. So either $f(v_{1}) \sim f(v_{3})$ or $f(v) \sim f(v_{2})$. Both cases are treated below.

e$.1.$ Suppose $f(v_{1}) \sim f(v_{3})$.  We choose a filling diagram $(D', f')$ for $\gamma$ such that $D'$ is triangulated with the same simplices like $D$ except for the triangles $\langle v,w,v_{1} \rangle$, $\langle w,v_{i},v_{i+1} \rangle, 1 \leq i \leq 2$  in $D$ which are replaced in
$D'$ by the triangles $\langle v,w,v_{3} \rangle$, $\langle v, v_{1}, v_{3} \rangle,$ $\langle v_{1}, v_{2}, v_{3} \rangle$.
We define $f'$ such that it coincides with $f$ on all simplices
which are common to $D$ and $D'$. We define $f'$ such that $f'(\langle v_{1}, v_{3} \rangle) = \langle f(v_{1}), f(v_{3}) \rangle$, $f'(\langle v, v_{3} \rangle) = \langle f(v), f(v_{3}) \rangle$, $f'(\langle v,w,v_{3} \rangle) = \langle f(v),f(w),f(v_{3}) \rangle$, $f'(\langle v,v_{1},v_{3} \rangle) = \langle f(v),f(v_{1}),f(v_{3}) \rangle$, $f'(\langle v_{1},v_{2},v_{3} \rangle) = \langle f(v_{1}),f(v_{2}),f(v_{3}) \rangle$.
Because $f$ is simplicial, so is $f'$. Then $(D', f')$ is indeed a filling diagram for $\gamma$.
 Because the discs $D$ and $D'$ have the same area, $D'$ has minimal area. Note that there is a full $4$-cycle $(v,v_{3},v_{4},v_{5})$ in $D'_{w}$ and there is a full $(k+1)$-cycle $(v_{1},v_{3},w,v_{5},...,v_{k+2})$, $4 \leq k \leq 6$ in $D'_{v}$. If $k=4$, then case $1$, $k=5$ (or case $2$, $k=4$) implies a contradiction.
If $k=5$, then case $1$, $k=6$ implies a contradiction.
Let $k=6$. Because $D'$ has minimal area, the map $f'$ is simplicial and nondegenerate. So there is a $4$-cycle $\eta_{1} = (f'(v),f'(v_{3}),f'(v_{4}),f'(v_{5}))$ in $X_{f'(w)}$ and there is a $7$-cycle $\eta_{2} = (f'(v_{1}),f'(v_{3}),f'(w),f'(v_{5}),f'(v_{6}),f'(v_{7}),f'(v_{8}))$ in $X_{f'(v)}$. Lemma \ref{3.2} implies that $\eta_{2}$ is not full. So at least one of the following holds $f'(v_{1}) \sim f'(w)$, $f'(v_{1}) \sim f'(v_{i})$, $5 \leq i \leq 7$. These cases are treated below.

$\bullet$ Suppose $f'(v_{1}) \sim f'(w)$. We choose a filling diagram $(D'',f'')$ for $\gamma$ such that $D''$ is triangulated with the same simplices like $D'$ except for the triangles $\langle v,v_{1},v_{3}\rangle$, $\langle v,w,v_{3} \rangle$ in $D'$ which are replaced in $D''$ by the triangles $\langle v,w,v_{1} \rangle$, $\langle w,v_{1},v_{3} \rangle$. We define $f''$ such that it coincides with $f'$ on all simplices which are common to $D'$ and $D''$. We define $f''$ such that $f''(\langle v_{1},w\rangle) = \langle f'(v_{1}),f'(w) \rangle$, $f''(\langle v,w,v_{1} \rangle) = \langle f'(v),f'(w),f'(v_{1}) \rangle$, $f''(\langle w,v_{1},v_{3} \rangle) = \langle f'(w),f'(v_{1}),f'(v_{3}) \rangle$. Because $f'$ is simplicial, $f''$ is also simplicial. So $(D'',f'')$ is indeed a filling diagram for $\gamma$. The area of $D''$ equals the area of $D'$ and therefore $D''$ has minimal area. Note that there is a full $5$-cycle $(v,v_{1},v_{3},v_{4},v_{5})$ in $D''_{w}$ and there is a full $6$-cycle $(v_{1},w,v_{5},v_{6},v_{7},v_{8})$ in $D''_{v}$. Then case $2$, $k=6$ implies a contradiction.

\begin{figure}[h]
    \begin{center}
       \includegraphics[height=4cm]{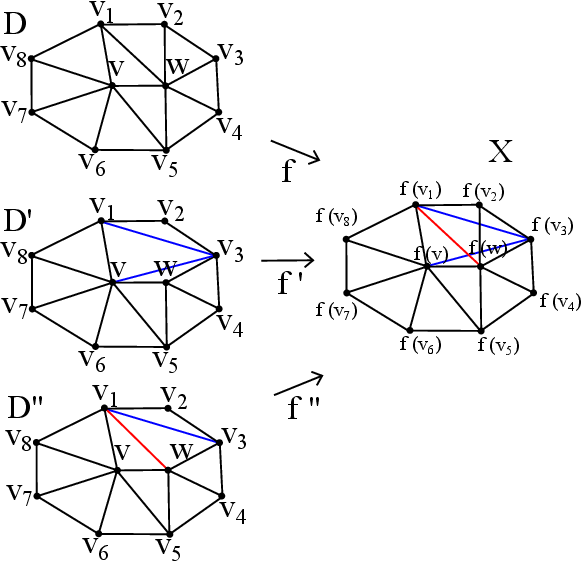}
      \caption{Case $3$, $k=6$, $f(v) \sim f(v_{3}), f(v_{1}) \sim f(v_{3}), f'(v_{1}) \sim f'(w)$}
    \end{center}
\end{figure}

$\bullet$ Suppose $f'(v_{1}) \sim f'(v_{5})$. By Lemma \ref{3.2}, the $5$-cycle $(f'(v_{1}),f'(v_{5}),$ $f'(v_{6}),$ $f'(v_{7}),f'(v_{8}))$ in $X_{f'(v)}$ is not full. So at least one of the following holds $f'(v_{1}) \sim f'(v_{i}), 6 \leq i \leq 7$. Assume w.l.o.g. $f'(v_{1}) \sim f'(v_{6})$. The other cases can be treated similarly. By Lemma \ref{3.2}, the $4$-cycle $(f'(v_{1}),f'(v_{6}),$ $f'(v_{7}),f'(v_{8}))$ in $X_{f'(v)}$ is not full.  Assume w.l.o.g. $f'(v_{1}) \sim f'(v_{7})$. We choose a filling diagram $(D'',f'')$ for $\gamma$ such that $D''$ is triangulated with the same simplices like $D'$ except for the triangles $\langle v,v_{i},v_{i+1}\rangle, 5 \leq i \leq 7$, $\langle v,v_{1},v_{8}\rangle$ in $D'$ which are replaced in $D''$ by the triangles $\langle v,v_{1},v_{5}\rangle$, $\langle v_{1},v_{i},v_{i+1}\rangle, 5 \leq i \leq 7$. We define $f''$ such that it coincides with $f'$ on all simplices which are common to $D'$and $D''$. We define $f''$ such that $f''(\langle v_{1},v_{i}\rangle) = \langle f'(v_{1}),f'(v_{i}) \rangle$, $f''(\langle v_{1},v_{i},v_{i+1} \rangle) = \langle f'(v_{1}),f'(v_{i}),f'(v_{i+1}) \rangle$, $5 \leq i \leq 7$, $f''(\langle v,v_{1},v_{5} \rangle) = \langle f'(v),f'(v_{1}),f'(v_{5}) \rangle$. Because $f'$ is simplicial, $f''$ is also simplicial. So $(D'',f'')$ is indeed a filling diagram for $\gamma$. The area of $D''$ equals the area of $D'$ and therefore $D''$ has minimal area. There is a full $4$-cycle $(v_{1},v_{3},w,v_{5})$ in $D''_{v}$, and a full $4$-cycle $(v,v_{3},v_{4},v_{5})$ in $D''_{w}$. Then case $1$, $k=4$ implies a contradiction.

\begin{figure}[h]
    \begin{center}
       \includegraphics[height=4cm]{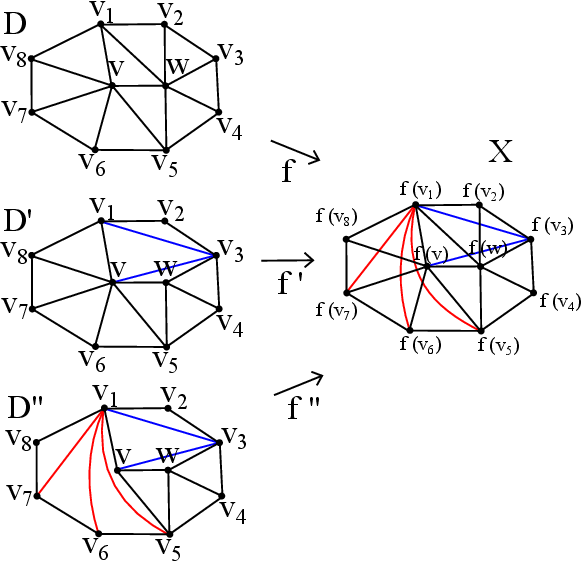}
      \caption{Case $3$, $k=6$, $f(v) \sim f(v_{3}), f(v_{1}) \sim f(v_{3}), f'(v_{1}) \sim f'(v_{5})$}
    \end{center}
\end{figure}

$\bullet$ Suppose $f'(v_{1}) \sim f'(v_{6})$. By Lemma \ref{3.2}, the $4$-cycle $(f'(v_{1}),f'(v_{6}),$ $f'(v_{7}),f'(v_{8}))$ in $X_{f'(v)}$ is not full. Assume w.l.o.g. $f'(v_{1}) \sim f'(v_{7})$. We choose a filling diagram $(D'',f'')$ for $\gamma$ such that $D''$ is triangulated with the same simplices like $D'$ except for the triangles $\langle v,v_{i},v_{i+1}\rangle, 6 \leq i \leq 7$, $\langle v,v_{1},v_{8}\rangle$ in $D'$ which are replaced in $D''$ by the triangles $\langle v,v_{1},v_{6}\rangle$, $\langle v_{1},v_{i},v_{i+1}\rangle, 6 \leq i \leq 7$. We define $f''$ such that it coincides with $f'$ on all simplices which are common to $D'$and $D''$. We define $f''$ such that $f''(\langle v_{1},v_{i}\rangle) = \langle f'(v_{1}),f'(v_{i}) \rangle$, $f''(\langle v_{1},v_{i},v_{i+1} \rangle) = \langle f'(v_{1}),f'(v_{i}),f'(v_{i+1}) \rangle$, $6 \leq i \leq 7$, $f''(\langle v,v_{1},v_{6} \rangle) = \langle f'(v),f'(v_{1}),f'(v_{6}) \rangle$. Because $f'$ is simplicial, $f''$ is also simplicial. So $(D'',f'')$ is indeed a filling diagram for $\gamma$. The area of $D''$ equals the area of $D'$ and therefore $D''$ has minimal area. There is a full $5$-cycle $(v_{1},v_{3},w,v_{5},v_{6})$ in $D''_{v}$, and a full $4$-cycle $(v,v_{3},v_{4},v_{5})$ in $D''_{w}$. Then case $1$, $k=5$ (or case $2$, $k=4$) implies a contradiction.

$\bullet$ Suppose $f'(v_{1}) \sim f'(v_{7})$. We choose a filling diagram $(D'',f'')$ for $\gamma$ such that $D''$ is triangulated with the same simplices like $D'$ except for the triangles $\langle v,v_{7},v_{8}\rangle$, $\langle v,v_{1},v_{8}\rangle$ in $D'$ are replaced in $D''$ by the triangles $\langle v_{1},v_{7},v_{8}\rangle$, $\langle v,v_{1},v_{7}\rangle$. We define $f''$ such that it coincides with $f'$ on all simplices which are common to $D'$ and $D''$. We define $f''$ such that $f''(\langle v_{1},v_{7}\rangle) = \langle f'(v_{1}),f'(v_{7}) \rangle$, $f''(\langle v_{1},v_{7},v_{8} \rangle) = \langle f'(v_{1}),f'(v_{7}),f'(v_{8}) \rangle$, \\ $f''(\langle v,v_{1},v_{7} \rangle) = \langle f'(v),f'(v_{1}),f'(v_{7}) \rangle$. Because $f'$ is simplicial, $f''$ is also simplicial. So $(D'',f'')$ is indeed a filling diagram for $\gamma$. The area of $D''$ equals the area of $D'$ and therefore $D''$ has minimal area. There is a full $6$-cycle $(v_{1},v_{3},w,v_{5},v_{6},v_{7})$ in $D''_{v}$, and a full $4$-cycle $(v,v_{3},v_{4},v_{5})$ in $D''_{w}$. Hence case $1$, $k=6$ implies a contradiction.

e$.2.$ Suppose $f(v) \sim f(v_{2})$. We choose a filling diagram $(D', f')$ for $\gamma$ such that $D'$ is triangulated with the same simplices like $D$ except for the triangles $\langle v,w,v_{1} \rangle, \langle w,v_{i},v_{i+1} \rangle, 1 \leq i \leq 2$  in $D$ which are replaced in
$D'$ by the triangles $\langle v, w, v_{3} \rangle$, $\langle v, v_{i}, v_{i+1} \rangle, 1 \leq i \leq 2$.
We define $f'$ such that it coincides with $f$ on all simplices
which are common to $D$ and $D'$. We define $f'$ such that $f'(\langle v, v_{i} \rangle) = \langle f(v), f(v_{i}) \rangle$, $2 \leq i \leq 3$, $f'(\langle v,v_{i},v_{i+1} \rangle) = \langle f(v),f(v_{i}),f(v_{i+1}) \rangle$, $1 \leq i \leq 2$, $f'(\langle v,w,v_{3} \rangle) = \langle f(v),f(w),f(v_{3}) \rangle$.
Because $f$ is simplicial, so is $f'$. Then $(D', f')$ is indeed a filling diagram for $\gamma$.
The discs $D$ and $D'$ have the same area. So $D'$ has minimal area. Note that there is a full $4$-cycle $(v,v_{3},v_{4},v_{5})$ in $D'_{w}$ and a full $(k+2)$-cycle $(v_{1},v_{2},v_{3},w,v_{5}, ..., v_{k+2})$ in $D'_{v}$, $4 \leq k \leq 6$.  If $k=4$, then case $1$, $k=6$ implies a contradiction. Because $f'$ is simplicial and nondegenerate, there is a $4$-cycle $\eta_{1} = (f'(v),f'(v_{3}),f'(v_{4}),f'(v_{5}))$ in $X_{f'(w)}$ and there is a $(k+2)$-cycle $\eta_{2} = (f'(v_{1}),f'(v_{2}),f'(v_{3}),f'(w),f'(v_{5}), ..., f'(v_{k+2}))$ in $X_{f'(v)}$, $5 \leq k \leq 6$.
 Lemma \ref{3.2} implies that $\eta_{2}$ is not full. So at least one of the following holds $f'(v_{1}) \sim f'(w)$, $f'(v_{1}) \sim f'(v_{3})$, $f'(v_{1}) \sim f'(v_{i}), 5 \leq i \leq k+1$, $5 \leq k \leq 6$. We treat some of these cases below. The other cases can be treated similarly.

$\bullet$ Let $k=5$. Suppose $f'(v_{1}) \sim f'(w)$. Note that there is a $5$-cycle $\eta_{3} = (f'(v_{1}),f'(w),f'(v_{5}),f'(v_{6}),f'(v_{7}))$ in $X_{f'(v)}$. Because $\eta_{1}$ is a $4$-cycle in $X_{f'(w)}$, Lemma \ref{3.2} implies that $\eta_{3}$ is not full. So at least one of the following holds $f'(v_{1}) \sim f'(v_{i})$, $5 \leq i \leq 6$. Assume w.l.o.g. $f'(v_{1}) \sim f'(v_{5})$. The other cases can be treated similarly. By Lemma \ref{3.2}, the $4$-cycle $(f'(v_{1}),f'(v_{5}),f'(v_{6}),$ $f'(v_{7}))$ in $X_{f'(v)}$ is not full. Assume w.l.o.g. $f'(v_{1}) \sim f'(v_{6})$. We choose a filling diagram $(D'',f'')$ for $\gamma$ such that $D''$ is triangulated with the same simplices like $D'$ except for the triangles $\langle v,w,v_{5} \rangle$, $\langle v,v_{i},v_{i+1}\rangle, 5 \leq i \leq 6$, $\langle v,v_{1},v_{7}\rangle$ in $D'$ which are replaced in $D''$ by the triangles $\langle v,w,v_{1}\rangle$, $\langle w,v_{1},v_{5} \rangle$, $\langle v_{1},v_{i},v_{i+1}\rangle, 5 \leq i \leq 6$. We define $f''$ such that it coincides with $f'$ on all simplices which are common to $D'$and $D''$. We define $f''$ such that $f''(\langle v_{1},w \rangle) = \langle f'(v_{1}),f'(w) \rangle$, $f''(\langle v_{1},v_{i}\rangle) = \langle f'(v_{1}),f'(v_{i}) \rangle$, $f''(\langle v_{1},v_{i},v_{i+1} \rangle) = \langle f'(v_{1}),f'(v_{i}),f'(v_{i+1}) \rangle$, $5 \leq i \leq 6$, $f''(\langle v,w,v_{1} \rangle) = \langle f'(v),f'(w),f'(v_{1}) \rangle$, $f''(\langle w,v_{1},v_{5} \rangle) =$ \\ $\langle f'(w),f'(v_{1}),f'(v_{5}) \rangle$. Because $f'$ is simplicial, $f''$ is also simplicial. So $(D'',f'')$ is indeed a filling diagram for $\gamma$. The area of $D''$ equals the area of $D'$ and therefore $D''$ has minimal area. Note that there is a full $4$-cycle $(v_{1},v_{2},v_{3},w)$ in $D''_{v}$ and a full $5$-cycle $(v_{1},v,v_{3},v_{4},v_{5})$ in $D''_{w}$. Then case $1$, $k=5$ implies a contradiction.

\begin{figure}[h]
    \begin{center}
       \includegraphics[height=4cm]{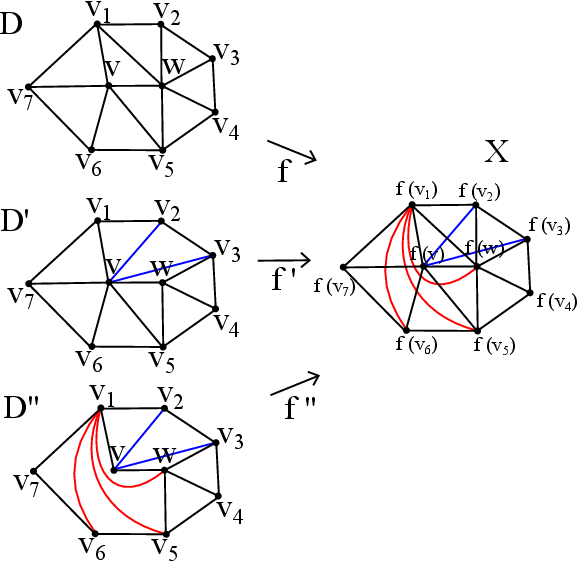}
      \caption{Case $3$, $k=5$, $f(v) \sim f(v_{3}), f(v) \sim f(v_{2}), f'(v_{1}) \sim f'(w)$}
    \end{center}
\end{figure}

$\bullet$ Let $k=5$. Suppose $f'(v_{1}) \sim f'(v_{3})$. We choose a filling diagram $(D'',f'')$ for $\gamma$ such that $D''$ is triangulated with the same simplices like $D'$ except for the triangles $\langle v,v_{i},v_{i+1} \rangle$, $1 \leq i \leq 2$ in $D'$ which are replaced in $D''$ by the triangles $\langle v,v_{1},v_{3}\rangle$, $\langle v_{1},v_{2},v_{3} \rangle$. We define $f''$ such that it coincides with $f'$ on all simplices which are common to $D'$and $D''$. We define $f''$ such that $f''(\langle v_{1},v_{3}\rangle) = \langle f'(v_{1}),f'(v_{3}) \rangle$, $f''(\langle v,v_{1},v_{3} \rangle) = \langle f'(v),f'(v_{1}),f'(v_{3}) \rangle$, $f''(\langle v_{1},v_{2},v_{3} \rangle) = \langle f'(v_{1}),f'(v_{2}),f'(v_{3}) \rangle$. Because $f'$ is simplicial, $f''$ is also simplicial. So $(D'',f'')$ is indeed a filling diagram for $\gamma$. The area of $D''$ equals the area of $D'$ and therefore $D''$ has minimal area. Note that there is a full $4$-cycle $(v,v_{3},v_{4},v_{5})$ in $D''_{w}$ and there is a full $6$-cycle $(v_{1},v_{3},w,v_{5},v_{6},v_{7})$ in $D''_{v}$. Then case $1$, $k=6$ implies a contradiction.

\begin{figure}[h]
    \begin{center}
       \includegraphics[height=4cm]{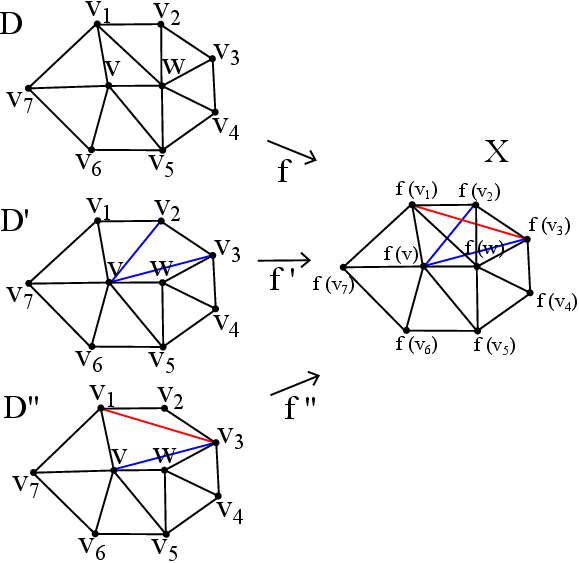}
      \caption{Case $3$, $k=5$, $f(v) \sim f(v_{3}), f(v) \sim f(v_{2}), f'(v_{1}) \sim f'(v_{3})$}
    \end{center}
\end{figure}

$\bullet$ Let $k=6$. Suppose $f'(v_{1}) \sim f'(w)$. Note that there is a $6$-cycle $\eta_{3} = (f'(v_{1}),f'(w),f'(v_{5}),f'(v_{6}),f'(v_{7}),f'(v_{8}))$ in $X_{f'(v)}$. Because $\eta_{1} = (f'(v),f'(v_{3}),f'(v_{4}),f'(v_{5}))$ is a $4$-cycle in $X_{f'(w)}$, Lemma \ref{3.2} implies that $\eta_{3}$ is not full. So at least one of the following holds $f'(v_{1}) \sim f'(v_{i})$, $5 \leq i \leq 7$. Assume w.l.o.g. $f'(v_{1}) \sim f'(v_{5})$. The other cases can be treated similarly. By Lemma \ref{3.2}, the $5$-cycle $(f'(v_{1}),f'(v_{5}),f'(v_{6}),$ $f'(v_{7}),f'(v_{8}))$ in $X_{f'(v)}$ is not full. So at least one of the following holds $f'(v_{1}) \sim f'(v_{i}), 6 \leq i \leq 7$. Assume w.l.o.g. $f'(v_{1}) \sim f'(v_{6})$. The other cases can be treated similarly. By Lemma \ref{3.2}, the $4$-cycle $(f'(v_{1}),f'(v_{6}),$ $f'(v_{7}),f'(v_{8}))$ in $X_{f'(v)}$ is not full.  Assume w.l.o.g. $f'(v_{1}) \sim f'(v_{7})$. We choose a filling diagram $(D'',f'')$ for $\gamma$ such that $D''$ is triangulated with the same simplices like $D'$ except for the triangles $\langle v,v_{i},v_{i+1} \rangle$, $5 \leq i \leq 7$, $\langle v,v_{1},v_{8}\rangle$, $\langle v,w,v_{5} \rangle$ in $D'$ which are replaced in $D''$ by the triangles $\langle v,w,v_{1}\rangle$, $\langle w,v_{1},v_{5}\rangle$, $\langle v_{1},v_{i},v_{i+1} \rangle$, $5 \leq i \leq 7$. We define $f''$ such that it coincides with $f'$ on all simplices which are common to $D'$ and $D''$. We define $f''$ such that $f''(\langle w,v_{1} \rangle) = \langle f'(w),f'(v_{1}) \rangle$, $f''(\langle v_{1},v_{i} \rangle) = \langle f'(v_{1}),f'(v_{i}) \rangle$, $5 \leq i \leq 7$, $f''(\langle v,w,v_{1} \rangle) = \langle f'(v),f'(w),f'(v_{1}) \rangle$, $f''(\langle w,v_{1},v_{5} \rangle) = \langle f'(w),f'(v_{1}),f'(v_{5}) \rangle$, \\ $f''(\langle v_{1},v_{i},v_{i+1} \rangle) = \langle f'(v_{1}),f'(v_{i}),f'(v_{i+1}) \rangle$, $5 \leq i \leq 7$. Because $f'$ is simplicial, $f''$ is also simplicial. So $(D'',f'')$ is indeed a filling diagram for $\gamma$. Since the area of $D''$ equals the area of $D'$, $D''$ has minimal area. Note that there is a full $4$-cycle $(w,v_{1},v_{2},v_{3})$ in $D''_{v}$ and a full $5$-cycle $(v_{1},v,v_{3},v_{4},v_{5})$ in $D''_{w}$. This implies a contradiction with case $1$, $k=5$ (or case $2$, $k=4$).

\begin{figure}[h]
    \begin{center}
       \includegraphics[height=4cm]{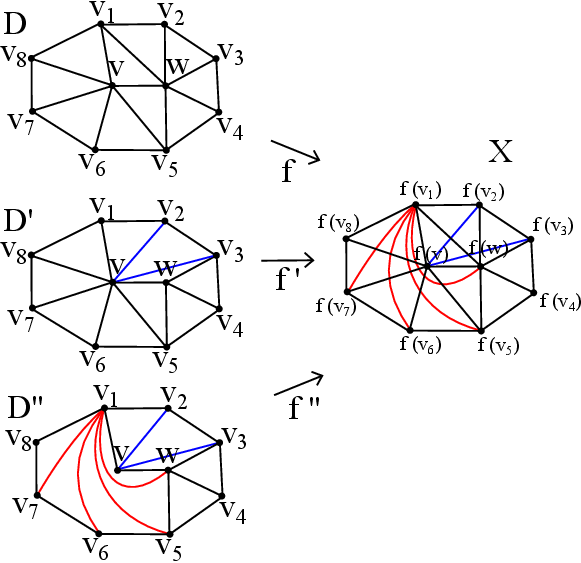}
      \caption{Case $3$, $k=6$, $f(v) \sim f(v_{3}), f(v) \sim f(v_{2}), f'(v_{1}) \sim f'(w)$}
    \end{center}
\end{figure}

$\bullet$ Let $k=6$. Suppose $f'(v_{1}) \sim f'(v_{3})$.  Note that there is a $7$-cycle $\eta_{4} = (f'(v_{1}),f'(v_{3}),f'(w),f'(v_{5}),f'(v_{6}),f'(v_{7}),f'(v_{8}))$ in $X_{f'(v)}$. Because $\eta_{1} = (f'(v),f'(v_{3}),f'(v_{4}),f'(v_{5}))$ is a $4$-cycle in $X_{f'(w)}$, Lemma \ref{3.2} implies that $\eta_{4}$ is not full. So at least one of the following holds $f'(v_{1}) \sim f'(w)$, $f'(v_{1}) \sim f'(v_{i}), 5 \leq i \leq 7$. These cases are treated below.

$\star$ Suppose $f'(v_{1}) \sim f'(w)$. We choose a filling diagram $(D'',f'')$ for $\gamma$ such that $D''$ is triangulated with the same simplices like $D'$ except for the triangles $\langle v,v_{i},v_{i+1}\rangle$, $1 \leq i \leq 2$, $\langle v,w,v_{3} \rangle$ in $D'$ which are replaced in $D''$ by the triangles $\langle v,w,v_{1}\rangle$, $\langle w,v_{1},v_{3} \rangle$, $\langle v_{1},v_{2},v_{3} \rangle$. We define $f''$ such that it coincides with $f'$ on all simplices which are common to $D'$ and $D''$. We define $f''$ such that $f''(\langle v_{1},w \rangle) = \langle f'(v_{1}),f'(w) \rangle$, $f''(\langle v_{1},v_{3} \rangle) = \langle f'(v_{1}),f'(v_{3}) \rangle$, $f''(\langle v,w,v_{1} \rangle) = \langle f'(v),f'(w),f'(v_{1}) \rangle$, $f''(\langle v_{1},v_{3},w \rangle) =$ \\ $\langle f'(v_{1}),f'(v_{3}),f'(w) \rangle$, $f''(\langle v_{1},v_{2},v_{3} \rangle) = \langle f'(v_{1}),f'(v_{2}),f'(v_{3}) \rangle$. Because $f'$ is simplicial, $f''$ is also simplicial. So $(D'',f'')$ is indeed a filling diagram for $\gamma$. The area of $D''$ equals the area of $D'$ and therefore $D''$ has minimal area. Note that there is a full $5$-cycle $(v,v_{1},v_{3},v_{4},v_{5})$ in $D''_{w}$ and a full $6$-cycle $(v_{1},w,v_{5},v_{6},v_{7},v_{8})$ in $D''_{v}$. Then case $2$, $k=6$ implies a contradiction.

\begin{figure}[h]
    \begin{center}
       \includegraphics[height=4cm]{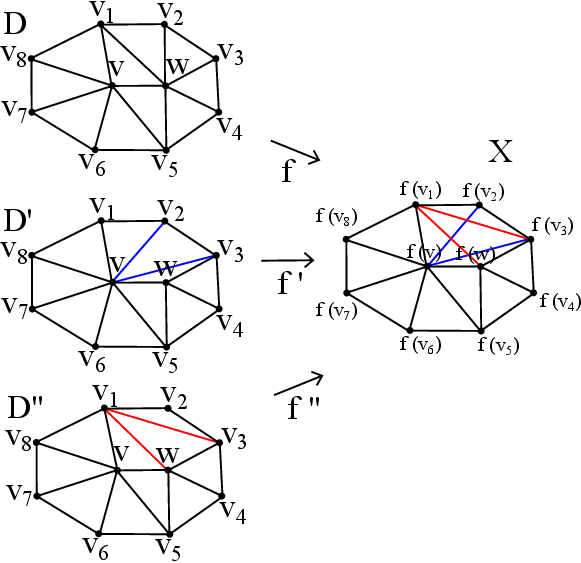}
      \caption{Case $3$, $k=6$, $f(v) \sim f(v_{3}), f(v) \sim f(v_{2}), f'(v_{1}) \sim f'(v_{3}), f'(v_{1}) \sim f'(w)$}
    \end{center}
\end{figure}

$\star$ Suppose $f'(v_{1}) \sim f'(v_{5})$. By Lemma \ref{3.2}, the $5$-cycle $(f'(v_{1}),f'(v_{5}),$ $f'(v_{6}),$ $f'(v_{7}),f'(v_{8}))$ in $X_{f'(v)}$ is not full. So at least one of the following holds $f'(v_{1}) \sim f'(v_{i}), 6 \leq i \leq 7$. Assume w.l.o.g. $f'(v_{1}) \sim f'(v_{6})$. The other cases can be treated similarly. By Lemma \ref{3.2}, the $4$-cycle $(f'(v_{1}),f'(v_{6}),$ $f'(v_{7}),f'(v_{8}))$ in $X_{f'(v)}$ is not full.  Assume w.l.o.g. $f'(v_{1}) \sim f'(v_{7})$. We choose a filling diagram $(D'',f'')$ for $\gamma$ such that $D''$ is triangulated with the same simplices like $D'$ except for the triangles $\langle v,v_{i},v_{i+1}\rangle$, $5 \leq i \leq 7$, $\langle v,v_{1},v_{8} \rangle$, $\langle v,v_{i},v_{i+1} \rangle$, $1 \leq i \leq 2$ in $D'$ which are replaced in $D''$ by the triangles $\langle v,v_{1},v_{5} \rangle$, $\langle v_{1},v_{i},v_{i+1} \rangle$, $5 \leq i \leq 7$, $\langle v,v_{1},v_{3} \rangle$, $\langle v_{1},v_{2},v_{3} \rangle$. We define $f''$ such that it coincides with $f'$ on all simplices which are common to $D'$ and $D''$. We define $f''$ such that $f''(\langle v_{1},v_{3}\rangle) = \langle f'(v_{1}),f'(v_{3}) \rangle$, $f''(\langle v_{1},v_{i}\rangle) = \langle f'(v_{1}),f'(v_{i}) \rangle$, $f''(\langle v_{1},v_{i},v_{i+1}\rangle) = \langle f'(v_{1}),f'(v_{i}),f'(v_{i+1}) \rangle$, $5 \leq i \leq 7$, $f''(\langle v,v_{1},v_{5} \rangle) = \langle f'(v),f'(v_{1}),f'(v_{5}) \rangle$, $f''(\langle v,v_{1},v_{3}\rangle) =$ \\ $\langle f'(v),f'(v_{1}),f'(v_{3}) \rangle$, $f''(\langle v_{1},v_{2},v_{3}\rangle) = \langle f'(v_{1}),f'(v_{2}),f'(v_{3}) \rangle$. Because $f'$ is simplicial, $f''$ is also simplicial. So $(D'',f'')$ is indeed a filling diagram for $\gamma$. The area of $D''$ equals the area of $D'$ and therefore $D''$ has minimal area. Note that there is a full $4$-cycle $(v,v_{3},v_{4},v_{5})$ in $D''_{w}$ and a full $4$-cycle $(v_{1},v_{3},w,v_{5})$ in $D''_{v}$. Then case $1$, $k=4$ implies a contradiction.

\begin{figure}[h]
    \begin{center}
       \includegraphics[height=4cm]{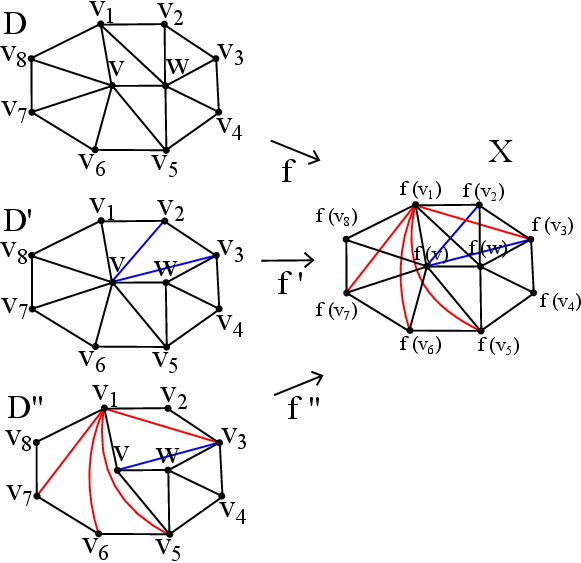}
      \caption{Case $3$, $k=6$, $f(v) \sim f(v_{3}), f(v) \sim f(v_{2}), f'(v_{1}) \sim f'(v_{3}), f'(v_{1}) \sim f'(v_{5})$}
    \end{center}
\end{figure}

$\star$ Suppose $f'(v_{1}) \sim f'(v_{6})$. By Lemma \ref{3.2}, the $4$-cycle $(f'(v_{1}),f'(v_{6}),$ $f'(v_{7}),f'(v_{8}))$ in $X_{f'(v)}$ is not full. Assume w.l.o.g. $f'(v_{1}) \sim f'(v_{7})$. We choose a filling diagram $(D'',f'')$ for $\gamma$ such that $D''$ is triangulated with the same simplices like $D'$ except for the triangles $\langle v,v_{i},v_{i+1}\rangle$, $6 \leq i \leq 7$, $\langle v,v_{1},v_{8} \rangle$, $\langle v,v_{i},v_{i+1} \rangle$, $1 \leq i \leq 2$ in $D'$ which are replaced in $D''$ by the triangles $\langle v,v_{1},v_{6}\rangle$, $\langle v_{1},v_{i},v_{i+1} \rangle$, $6 \leq i \leq 7$, $\langle v,v_{1},v_{3} \rangle$, $\langle v_{1},v_{2},v_{3} \rangle$. We define $f''$ such that it coincides with $f'$ on all simplices which are common to $D'$ and $D''$. We define $f''$ such that $f''(\langle v_{1},v_{3}\rangle) = \langle f'(v_{1}),f'(v_{3}) \rangle$, $f''(\langle v_{1},v_{i}\rangle) = \langle f'(v_{1}),f'(v_{i}) \rangle$, $f''(\langle v_{1},v_{i},v_{i+1}\rangle) = \langle f'(v_{1}),f'(v_{i}),f'(v_{i+1}) \rangle$, $6 \leq i \leq 7$, $f''(\langle v,v_{1},v_{6} \rangle) = \langle f'(v),f'(v_{1}),f'(v_{6}) \rangle$, $f''(\langle v,v_{1},v_{3} \rangle) =$ \\ $\langle f'(v),f'(v_{1}),f'(v_{3}) \rangle$, $f''(\langle v_{1},v_{2},v_{3} \rangle) = \langle f'(v_{1}),f'(v_{2}),f'(v_{3}) \rangle$. Because $f'$ is simplicial, $f''$ is also simplicial. So $(D'',f'')$ is indeed a filling diagram for $\gamma$. The area of $D''$ equals the area of $D'$ and therefore $D''$ has minimal area. Note that there is a full $4$-cycle $(v,v_{3},v_{4},v_{5})$ in $D''_{w}$ and a full $5$-cycle $(v_{1},v_{3},w,v_{5},v_{6})$ in $D''_{v}$. Then case $1$, $k=5$ implies a contradiction.

\begin{figure}[h]
    \begin{center}
       \includegraphics[height=4cm]{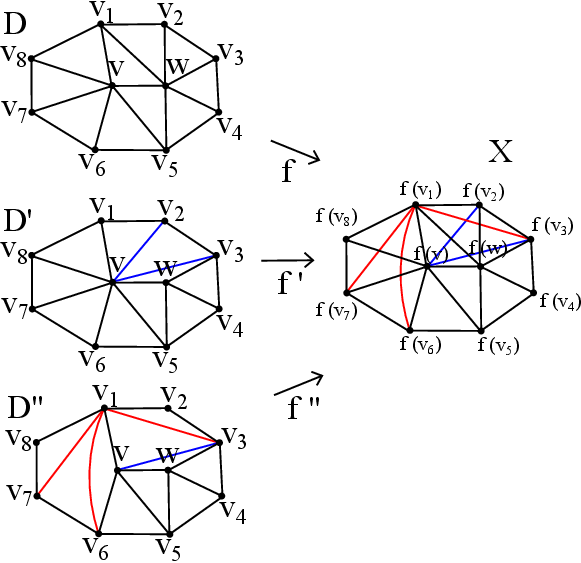}
      \caption{Case $3$, $k=6$, $f(v) \sim f(v_{3}), f(v) \sim f(v_{2}), f'(v_{1}) \sim f'(v_{3}), f'(v_{1}) \sim f'(v_{6})$}
    \end{center}
\end{figure}

$\star$ Suppose $f'(v_{1}) \sim f'(v_{7})$. We choose a filling diagram $(D'',f'')$ for $\gamma$ such that $D''$ is triangulated with the same simplices like $D'$ except for the triangles $\langle v,v_{7},v_{8}\rangle$, $\langle v,v_{1},v_{8} \rangle$, $\langle v,v_{i},v_{i+1} \rangle$, $1 \leq i \leq 2$ in $D'$ which are replaced in $D''$ by the triangles $\langle v,v_{1},v_{7}\rangle$, $\langle v_{1},v_{7},v_{8} \rangle$, $\langle v,v_{1},v_{3} \rangle$, $\langle v_{1},v_{2},v_{3} \rangle$. We define $f''$ such that it coincides with $f'$ on all simplices which are common to $D'$ and $D''$. We define $f''$ such that $f''(\langle v_{1},v_{7}\rangle) = \langle f'(v_{1}),f'(v_{7}) \rangle$, $f''(\langle v_{1},v_{3}\rangle) = \langle f'(v_{1}),f'(v_{3}) \rangle$, $f''(\langle v,v_{1},v_{7} \rangle) = \langle f'(v),f'(v_{1}),f'(v_{7}) \rangle$, \\ $f''(\langle v_{1},v_{7},v_{8} \rangle) = \langle f'(v_{1}),f'(v_{7}),f'(v_{8}) \rangle$, $f''(\langle v,v_{1},v_{3} \rangle) =$ \\ $\langle f'(v),f'(v_{1}),f'(v_{3}) \rangle$, $f''(\langle v_{1},v_{2},v_{3} \rangle) = \langle f'(v_{1}),f'(v_{2}),f'(v_{3}) \rangle$. Because $f'$ is simplicial, $f''$ is also simplicial. So $(D'',f'')$ is indeed a filling diagram for $\gamma$. The area of $D''$ equals the area of $D'$ and therefore $D''$ has minimal area. Note that there is a full $4$-cycle $(v,v_{3},v_{4},v_{5})$ in $D''_{w}$ and a full $6$-cycle $(v_{1},v_{3},w,v_{5},v_{6},v_{7})$ in $D''_{v}$. Then case $1$, $k=6$ implies a contradiction.

\begin{figure}[h]
    \begin{center}
       \includegraphics[height=4cm]{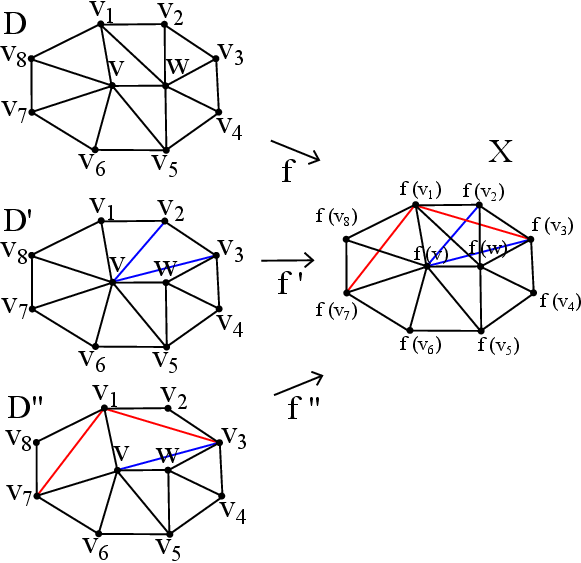}
      \caption{Case $3$, $k=6$, $f(v) \sim f(v_{3}), f(v) \sim f(v_{2}), f'(v_{1}) \sim f'(v_{3}), f'(v_{1}) \sim f'(v_{7})$}
    \end{center}
\end{figure}

So $w$ is $7$-large.

\item Case $4$: $k'' = 7$.

Let $D_{v} = (v_{1}, w, v_{6}, ..., v_{k+3})$, $4 \leq k \leq 5$ and let $D_{w} = (v_{1}, v_{2}, v_{3}, v_{4}, v_{5},$ $ v_{6}, v)$.
Lemma \ref{3.2} implies that the $7$-cycle
$(f(v_{1}), f(v_{2}), f(v_{3}), f(v_{4}), f(v_{5}),$ $ f(v_{6}), $ $f(v))$ in $X_{f(w)}$ is not full. We distinguish the cases $f(v_{1}) \sim f(v_{i})$, $3 \leq i \leq 6$, $f(v) \sim f(v_{i})$, $2 \leq i \leq 3$ which are treated below.

a. Suppose $f(v_{1}) \sim f(v_{3})$. We choose a filling diagram $(D', f')$ for $\gamma$ such that $D'$ is triangulated with the same simplices like $D$ except for the triangles $\langle w, v_{i}, v_{i+1} \rangle, 1 \leq i \leq 2$  in $D$ which are replaced in
$D'$ by the triangles $\langle w, v_{1}, v_{3} \rangle,$ $\langle v_{1}, v_{2}, v_{3} \rangle$.
We define $f'$ such that it coincides with $f$ on all simplices
which are common to $D$ and $D'$. We define $f'$ such that $f'(\langle v_{1}, v_{3} \rangle) = \langle f(v_{1}), f(v_{3}) \rangle$, $f'(\langle w,v_{1},v_{3} \rangle) = \langle f(w),f(v_{1}),f(v_{3}) \rangle$, $f'(\langle v_{1},v_{2},v_{3} \rangle) = $ $\langle f(v_{1}),$ $f(v_{2}),f(v_{3}) \rangle$.
Because $f$ is simplicial, so is $f'$. Then $(D', f')$ is indeed a filling diagram for $\gamma$.
The discs $D$ and $D'$ have the same area. So $D'$ has minimal area. Note that there is a full $6$-cycle $(v_{1},v_{3},v_{4},v_{5},v_{6},v)$ in $D'_{w}$ and there is a full $k$-cycle $(v_{1},w,v_{6}, ... ,v_{k+3})$, $4 \leq k \leq 5$ in $D'_{v}$. If $k=4$, then case $1$, $k=6$ implies a contradiction.  If $k=5$, we obtain contradiction with case $2$, $k=6$.

b. Suppose $f(v_{1}) \sim f(v_{4})$.
Lemma \ref{3.2} implies that the $4$-cycle $(f(v_{1}),$ $f(v_{2}),$ $ f(v_{3}), f(v_{4}))$ in $X_{f(w)}$ is not full. Assume w.l.og. $f(v_{1}) \sim f(v_{3})$.
 We choose a filling diagram $(D', f')$ for $\gamma$ such that $D'$ is triangulated with the same simplices like $D$
except for the triangles $\langle w, v_{i}, v_{i+1} \rangle, 1 \leq i \leq 3$ in $D$ which are replaced in $D'$ by the triangles
$\langle w, v_{1}, v_{4} \rangle, \langle v_{1}, v_{i}, v_{i+1} \rangle, 2 \leq i \leq 3$.
 We define $f'$ such that it coincides with $f$ on all simplices which are common to $D$ and $D'$. We define $f'$ such that
$f'(\langle v_{1}, v_{i} \rangle) = \langle f(v_{1}), f(v_{i}) \rangle$, $3 \leq i \leq 4$, $f'(\langle v_{1},v_{i},v_{i+1} \rangle) = \langle f(v_{1}),f(v_{i}),f(v_{i+1}) \rangle$, $2 \leq i \leq 3$, $f'(\langle w,v_{1},v_{4} \rangle) = \langle f(w),f(v_{1}),f(v_{4}) \rangle$.
Since $f$ is simplicial, so is $f'$. Then $(D', f')$ is indeed a filling diagram for $\gamma$. The discs $D$ and $D'$ have the same area.
So $D'$ has minimal area.
There is a full $5$-cycle $(v_{1}, v_{4}, v_{5}, v_{6}, v)$ in $D'_{w}$ and there is a full $k$-cycle $(v_{1}, w, v_{6}, ..., v_{k+3}), 4 \leq k \leq 5$ in $D'_{v}$.
If $k=4$, then
case $2$, $k=4$ (or case $1$, $k=5$) implies a contradiction.
If $k=5$, then
case $2$, $k=5$ implies a contradiction.

\begin{figure}[h]
    \begin{center}
       \includegraphics[height=4cm]{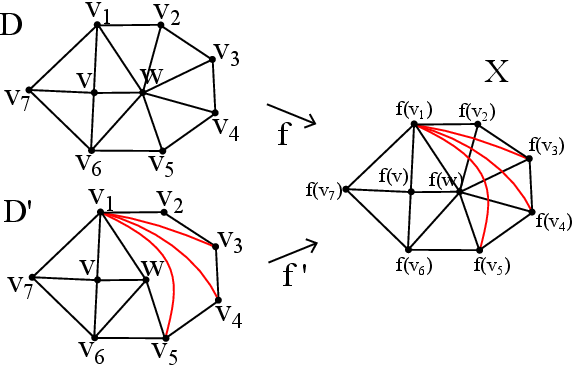}
      \caption{Case $4$, $k = 4$, $f(v_{1}) \sim f(v_{5})$}
    \end{center}
\end{figure}

\begin{figure}[h]
    \begin{center}
       \includegraphics[height=4cm]{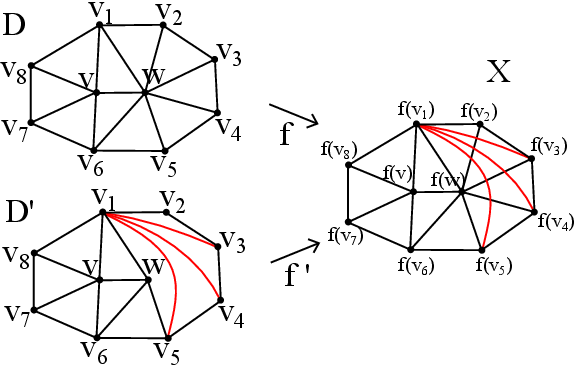}
      \caption{Case $4$, $k = 5$, $f(v_{1}) \sim f(v_{5})$}
    \end{center}
\end{figure}

c. Suppose $f(v_{1}) \sim f(v_{5})$.
Lemma \ref{3.2} implies that the $5$-cycle $(f(v_{1}),$ $f(v_{2}),$ $ f(v_{3}), f(v_{4}), f(v_{5}))$ in $X_{f(w)}$ is not full. So at least one of the following holds $f(v_{1}) \sim f(v_{i})$, $3 \leq i \leq 4$. Assume w.l.o.g. $f(v_{1}) \sim f(v_{4})$. The other cases can be treated similarly. Lemma \ref{3.2} implies that the $4$-cycle $(f(v_{1}), f(v_{2}), f(v_{3}), f(v_{4}))$ in $X_{f(w)}$ is not full. Assume w.l.o.g. $f(v_{1}) \sim f(v_{3})$.
 We choose a filling diagram $(D', f')$ for $\gamma$ such that $D'$ is triangulated with the same simplices like $D$
except for the triangles $\langle w, v_{i}, v_{i+1} \rangle, 1 \leq i \leq 4$ in $D$ which are replaced in $D'$ by the triangles
$\langle w, v_{1}, v_{5} \rangle, \langle v_{1}, v_{i}, v_{i+1} \rangle, 2 \leq i \leq 4$.
 We define $f'$ such that it coincides with $f$ on all simplices which are common to $D$ and $D'$. We define $f'$ such that $f'(\langle v_{1}, v_{i} \rangle) = \langle f(v_{1}), f(v_{i}) \rangle$, $3 \leq i \leq 5$, $f'(\langle v_{1},v_{i},v_{i+1} \rangle) = \langle f(v_{1}),f(v_{i}),f(v_{i+1}) \rangle, 2 \leq i \leq 4$, $f'(\langle w,v_{1},v_{5} \rangle) = \langle f(w),f(v_{1}),f(v_{5}) \rangle$.
Since $f$ is simplicial, so is $f'$. Then $(D', f')$ is indeed a filling diagram for $\gamma$. The discs $D$ and $D'$ have the same area.
So $D'$ has minimal area.
There is a full $4$-cycle $(v_{1}, v_{5}, v_{6}, v)$ in $D'_{w}$ and there is a full $k$-cycle $(v_{1}, w, v_{6}, ..., v_{k+3}), 4 \leq k \leq 5$ in $D'_{v}$.
If $k=4$, then
case $1$, $k=4$ implies a contradiction.
If $k=5$, then case $1$, $k=5$
(or case $2$, $k=4$) implies a contradiction.

d. Suppose $f(v_{1}) \sim f(v_{6})$.

$\bullet$ Let $k=4$. We consider a filling diagram $(D',f')$ for $\gamma$ such that $D'$ is triangulated with the same simplices like $D$ except for the triangles $\langle v,w,v_{1} \rangle$, $\langle v,w,v_{6} \rangle$, $\langle v,v_{6},v_{7} \rangle$, $\langle v,v_{1},v_{7} \rangle$ in $D$ which are replaced in $D'$ by the triangles $\langle w,v_{1},v_{6} \rangle$, $\langle v_{1},v_{6},v_{7} \rangle$. We define $f'$ such that it coincides with $f$ on all simplices which are common to $D$ and $D'$. We define $f'$ such that $f'(\langle v_{1},v_{6} \rangle) = \langle f(v_{1}),f(v_{6}) \rangle$, $f'(\langle w,v_{1},v_{6} \rangle) = \langle f(w),f(v_{1}),f(v_{6}) \rangle$, $f'(\langle v_{1},v_{6},$ $v_{7} \rangle) = \langle f(v_{1}),f(v_{6}),f(v_{7}) \rangle$. Because $f$ is simplicial, $f'$ is also simplicial. So $(D',f')$ is a filling diagram for $\gamma$. Note that the area of $D'$ is less than the area of $D$. Because $D$ has minimal area, we have reached a contradiction.

$\bullet$ Let $k=5$. Note that there is a $4$-cycle $\alpha = (f(v_{1}),f(v_{6}),f(v_{7}),f(v_{8}))$ in $X_{f(v)}$ and there is a $6$-cycle $(f(v_{1}), f(v_{2}), f(v_{3}), f(v_{4}), f(v_{5}), f(v_{6}))$ in $X_{f(w)}$. Lemma \ref{3.2} implies that $\alpha$ is not full. Assume w.l.o.g. $f(v_{1}) \sim f(v_{7})$. We consider a filling diagram for $\gamma$ such that $D'$ is triangulated with the same simplices like $D$ except for the triangles $\langle v,v_{i},v_{i+1} \rangle$, $6 \leq i \leq 7$, $\langle v,v_{1},v_{8} \rangle$ in $D$ which are replaced in $D'$ by the triangles $\langle w,v_{1},v_{6} \rangle$, $\langle v_{1},v_{i},v_{i+1} \rangle$, $6 \leq i \leq 7$. We define $f'$ such that $f'(\langle v_{1},v_{i} \rangle) = \langle f(v_{1}),f(v_{i}) \rangle$, $6 \leq i \leq 7$, $f'(\langle w,v_{1},v_{6} \rangle) = \langle f(w),f(v_{1}),f(v_{6}) \rangle$, $f'(v_{1},v_{i},v_{i+1}) =$ \\ $\langle f(v_{1}),f(v_{i}),f(v_{i+1}) \rangle$, $6 \leq i \leq 7$. Because $f$ is simplicial, $f'$ is also simplicial. So $(D',f')$ is a filling diagram for $\gamma$. Note that the area of $D'$ is less than the area of $D$. Because $D$ has minimal area, we have reached a contradiction.

e. Suppose $f(v) \sim f(v_{2})$. We choose a filling diagram $(D', f')$ for $\gamma$ such that $D'$ is triangulated with the same simplices like $D$ except for the triangles $\langle w,v,v_{1} \rangle, \langle w,v_{1},v_{2} \rangle$  in $D$ which are replaced in
$D'$ by the triangles $\langle v, w, v_{2} \rangle$, $\langle v, v_{1}, v_{2} \rangle$.
We define $f'$ such that it coincides with $f$ on all simplices
which are common to $D$ and $D'$. We define $f'$ such that $f'(\langle v, v_{2} \rangle) = \langle f(v), f(v_{2}) \rangle$, $f'(\langle v,w,v_{2} \rangle) = \langle f(v),f(w),f(v_{2}) \rangle$, $f'(\langle v,v_{1},v_{2} \rangle) = $ $\langle f(v),$ $f(v_{1}),f(v_{2}) \rangle$.
Because $f$ is simplicial, so is $f'$. Then $(D', f')$ is indeed a filling diagram for $\gamma$.
 The discs $D$ and $D'$ have the same area. So $D'$ has minimal area. Note that there is a full $6$-cycle $(v,v_{2},v_{3},v_{4},v_{5},v_{6})$ in $D'_{w}$ and there is a full $(k+1)$-cycle $(v_{1},v_{2},w,v_{6},...,v_{k+3})$ in $D'_{v}$, $4 \leq k \leq 5$. If $k=4$, then case $2$, $k=6$ (case $3$, $k=5$) implies a contradiction.  If $k=5$, then we obtain contradiction with case $3$, $k=6$.

f. Suppose $f(v) \sim f(v_{3})$. Note that there is a $4$-cycle $\alpha = (f(v),f(v_{1}),$ $f(v_{2}),f(v_{3}))$ in $X_{f(w)}$ and there is a $k$-cycle $(f(v_{1}),f(w),f(v_{6}), ... ,f(v_{k+3}))$ in $X_{f(v)}$, $4 \leq k \leq 5$.
Lemma \ref{3.2} implies that $\alpha$ is not full. So either $f(v_{1}) \sim f(v_{3})$ or $f(v) \sim f(v_{2})$. Both cases are treated below.

f.$1.$ Suppose $f(v_{1}) \sim f(v_{3})$.  We choose a filling diagram $(D', f')$ for $\gamma$ such that $D'$ is triangulated with the same simplices like $D$ except for the triangles $\langle v,w,v_{1} \rangle$, $\langle w,v_{i},v_{i+1} \rangle, 1 \leq i \leq 2$  in $D$ which are replaced in
$D'$ by the triangles $\langle v,w,v_{3} \rangle$, $\langle v, v_{1}, v_{3} \rangle,$ $\langle v_{1}, v_{2}, v_{3} \rangle$.
We define $f'$ such that it coincides with $f$ on all simplices
which are common to $D$ and $D'$. We define $f'$ such that $f'(\langle v_{1}, v_{3} \rangle) = \langle f(v_{1}), f(v_{3}) \rangle$, $f'(\langle v, v_{3} \rangle) = \langle f(v), f(v_{3}) \rangle$, $f'(\langle v,w,v_{3} \rangle) = \langle f(v),f(w),f(v_{3}) \rangle$, $f'(\langle v,v_{1},v_{3} \rangle) = \langle f(v),f(v_{1}),f(v_{3}) \rangle$, $f'(\langle v_{1},v_{2},v_{3} \rangle) = \langle f(v_{1}),f(v_{2}),f(v_{3}) \rangle$.
Because $f$ is simplicial, so is $f'$. Then $(D', f')$ is indeed a filling diagram for $\gamma$.
The discs $D$ and $D'$ have the same area. So $D'$ has minimal area. Note that there is a full $(k+1)$-cycle $(v_{1},v_{3},w,v_{6},...,v_{k+3})$, $4 \leq k \leq 5$ in $D'_{v}$ and there is a full $5$-cycle $(v,v_{3},v_{4},v_{5},v_{6})$ in $D'_{w}$. If $k=4$, then case $2$, $k=5$ implies a contradiction.
If $k=5$, then case $2$, $k=6$ (or case $3$, $k=5$) implies a contradiction.

 f.$2.$ Suppose $f(v) \sim f(v_{2})$. We choose a filling diagram $(D', f')$ for $\gamma$ such that $D'$ is triangulated with the same simplices like $D$ except for the triangles $\langle w,v,v_{1} \rangle, \langle w,v_{i},v_{i+1} \rangle, 1 \leq i \leq 2$  in $D$ which are replaced in
$D'$ by the triangles $\langle v, w, v_{3} \rangle$, $\langle v, v_{i}, v_{i+1} \rangle, 1 \leq i \leq 2$.
We define $f'$ such that it coincides with $f$ on all simplices
which are common to $D$ and $D'$. We define $f'$ such that $f'(\langle v, v_{i} \rangle) = \langle f(v), f(v_{i}) \rangle$, $2 \leq i \leq 3$, $f'(\langle v,v_{i},v_{i+1} \rangle) = \langle f(v),f(v_{i}),f(v_{i+1}) \rangle$, $1 \leq i \leq 2$, $f'(\langle v,w,v_{3} \rangle) = \langle f(v),f(w),f(v_{3}) \rangle$.
Because $f$ is simplicial, so is $f'$. Then $(D', f')$ is indeed a filling diagram for $\gamma$.
The discs $D$ and $D'$ have the same area. So $D'$ has minimal area. Note that there is a full $5$-cycle $(v,v_{3},v_{4},v_{5},v_{6})$ in $D'_{w}$ and a full $(k+2)$-cycle $(v_{1},v_{2},v_{3},w,v_{6},,...,v_{k+3})$, $4 \leq k \leq 6$ in $D'_{v}$. If $k=4$ then case $2$, $k=6$ (or case $3$, $k=5$) implies a contradiction.
Let $k=5$.
 Note that there is a $5$-cycle $\eta_{1} = (f'(v),f'(v_{3}),f'(v_{4}),f'(v_{5}),f'(v_{6}))$ in $X_{f'(w)}$ and there is a $7$-cycle $\eta_{2} = (f'(v_{1}),f'(v_{2}),$ $f'(v_{3}),f'(w),f'(v_{6}),f'(v_{7}),f'(v_{8}))$ in $X_{f'(v)}$. Lemma \ref{3.2} implies that $\eta_{2}$ is not full. So at least one of the following holds $f'(v_{1}) \sim f'(v_{3})$, $f'(v_{1}) \sim f'(w)$, $f'(v_{1}) \sim f'(v_{i})$, $6 \leq i \leq 7$. We treat below some of these cases. The other cases can be treated similarly.


\begin{figure}[h]
    \begin{center}
       \includegraphics[height=4cm]{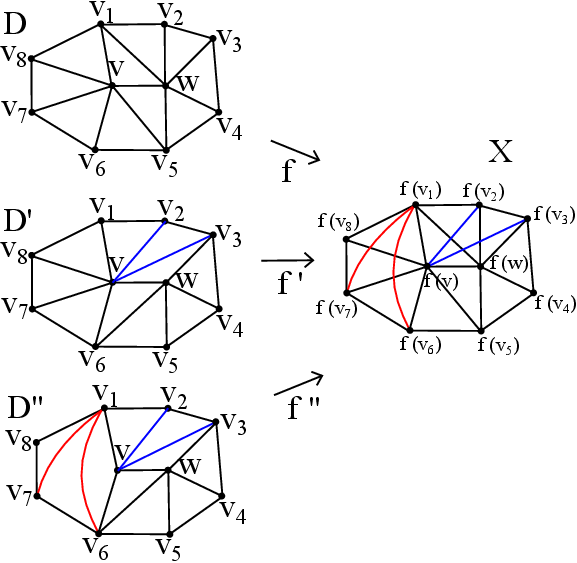}
      \caption{Case $4$, $k=5$, $f(v) \sim f(v_{3})$, $f(v) \sim f(v_{2})$, $f'(v_{1}) \sim f'(v_{6})$}
    \end{center}
\end{figure}

$\bullet$ Suppose $f'(v_{1}) \sim f'(v_{6})$. Note that there is a $4$-cycle $\alpha = (f'(v_{1}),f'(v_{6}),$ $f'(v_{7}),f'(v_{8}))$ in $X_{f'(v)}$ and there is a $5$-cycle $(f'(v),f'(v_{3}),$ $f'(v_{4}),f'(v_{5}),$ $f'(v_{6}))$ in $X_{f'(w)}$. By Lemma \ref{3.2}, $\alpha$ is not full. Assume w.l.o.g. $f'(v_{1}) \sim f'(v_{7})$. We choose a filling diagram $(D'',f'')$ for $\gamma$ such that $D''$ is triangulated with the same simplices like $D'$ except for the triangles $\langle v,v_{i},v_{i+1}\rangle, 6 \leq i \leq 7$, $\langle v,v_{1},v_{8}\rangle$ in $D'$ are replaced in $D''$ by the triangles $\langle v,v_{1},v_{6}\rangle$, $\langle v_{1},v_{i},v_{i+1}\rangle, 6 \leq i \leq 7$. We define $f''$ such that it coincides with $f'$ on all simplices which are common to $D'$and $D''$. We define $f''$ such that $f''(\langle v_{1},v_{i}\rangle) = \langle f'(v_{1}),f'(v_{i}) \rangle$, $f''(\langle v_{1},v_{i},v_{i+1} \rangle) = \langle f'(v_{1}),f'(v_{i}),f'(v_{i+1}) \rangle$, $6 \leq i \leq 7$, $f''(\langle v,v_{1},v_{6} \rangle) = \langle f'(v),f'(v_{1}),f'(v_{6}) \rangle$. Because $f'$ is simplicial, $f''$ is also simplicial. So $(D'',f'')$ is indeed a filling diagram for $\gamma$. The area of $D''$ equals the area of $D'$ and therefore $D''$ has minimal area. There is a full $5$-cycle $(v_{1},v_{2},v_{3},w,v_{6})$ in $D''_{v}$, and a full $5$-cycle $(v,v_{3},v_{4},v_{5},v_{6})$ in $D''_{w}$. Then case $2$, $k=5$ implies a contradiction.

$\bullet$ Suppose $f'(v_{1}) \sim f'(v_{7})$. We choose a filling diagram $(D'',f'')$ for $\gamma$ such that $D''$ is triangulated with the same simplices like $D'$ except for the triangles $\langle v,v_{7},v_{8}\rangle$, $\langle v,v_{1},v_{8}\rangle$ in $D'$ which are replaced in $D''$ by the triangles $\langle v,v_{1},v_{7}\rangle$, $\langle v_{1},v_{7},v_{8}\rangle$. We define $f''$ such that it coincides with $f'$ on all simplices which are common to $D'$ and $D''$. We define $f''$ such that $f''(\langle v_{1},v_{7}\rangle) = \langle f'(v_{1}),f'(v_{7}) \rangle$, $f''(\langle v,v_{1},v_{7} \rangle) = \langle f'(v),f'(v_{1}),f'(v_{7}) \rangle$, $f''(\langle v_{1},v_{7},v_{8} \rangle) = \langle f'(v_{1}),f'(v_{7}),f'(v_{8}) \rangle$. Because $f'$ is simplicial, $f''$ is also simplicial. So $(D'',f'')$ is indeed a filling diagram for $\gamma$. The area of $D''$ equals the area of $D'$ and therefore $D''$ has minimal area. There is a full $6$-cycle $(v_{1},v_{2},v_{3},w,v_{6},v_{7})$ in $D''_{v}$, and a full $5$-cycle $(v,v_{3},v_{4},v_{5},v_{6})$ in $D''_{w}$. Hence case $2$, $k=6$ (or case $3$, $k=5$) implies a contradiction.



Thus $w$ is $8$-large.

\item Case $5$: $k'' = 8$.

Let $D_{v} = (v_{1}, w, v_{7}, v_{8})$ and
let $D_{w} = (v_{1}, v_{2}, v_{3}, v_{4}, v_{5}, v_{6}, v_{7}, v)$.
Lemma \ref{3.2} implies that the $8$-cycle $(f(v_{1}), f(v_{2}), f(v_{3}), f(v_{4}),$ $ f(v_{5}),$ $ f(v_{6}), f(v_{7}),$ $f(v))$ in $X_{f(w)}$ is not full. We distinguish the cases $f(v_{1}) \sim f(v_{i}), 3 \leq i \leq 7$, $f(v) \sim f(v_{i})$, $2 \leq i \leq 4$.

a. Suppose $f(v_{1}) \sim f(v_{3})$. We choose a filling diagram $(D', f')$ for $\gamma$ such that $D'$ is triangulated with the same simplices like $D$ except for the triangles $\langle w, v_{i}, v_{i+1} \rangle, 1 \leq i \leq 2$ in $D$ which are replaced in
$D'$ by the triangles $\langle w, v_{1}, v_{3} \rangle$, $\langle v_{1}, v_{2}, v_{3} \rangle$.
We define $f'$ such that it coincides with $f$ on all simplices
which are common to $D$ and $D'$. We define $f'$ such that $f'(\langle v_{1}, v_{3} \rangle) = \langle f(v_{1}), f(v_{3}) \rangle$,
$f'(\langle w,v_{1},v_{3} \rangle) = \langle f(w),f(v_{1}),f(v_{3}) \rangle$, $f'(\langle v_{1},v_{2},v_{3} \rangle) = \langle f(v_{1}),$ $f(v_{2}),f(v_{3}) \rangle$. Because $f$ is simplicial, so is $f'$. Then $(D', f')$ is indeed a filling diagram for $\gamma$.
Note that the discs $D$ and $D'$ have the same area. So $D'$ has minimal area. There is a full $4$-cycle $(v_{1},w,v_{7},v_{8})$ in $D'_{v}$, and a full $7$-cycle  $(v,v_{1},v_{3},v_{4},v_{5},v_{6},v_{7})$ in $D'_{w}$. Hence case $4$, $k=4$ implies a contradiction.

b. Suppose $f(v_{1}) \sim f(v_{4})$. Note that there is a $4$-cycle $\alpha = (f(v_{1}),$ $f(v_{2}),f(v_{3}),f(v_{4}))$ in $X_{f(w)}$ and there is a $4$-cycle $(f(v_{1}),$ $f(w),f(v_{7}),f(v_{8}))$ in $X_{f(v)}$. According to Lemma \ref{3.2}, $\alpha$ is not full. Assume w.l.o.g. $f(v_{1}) \sim f(v_{3})$. We choose a filling diagram $(D', f')$ for $\gamma$ such that $D'$ is triangulated with the same simplices like $D$ except for the triangles $\langle w, v_{i}, v_{i+1} \rangle, 1 \leq i \leq 3$ in $D$ which are replaced in
$D'$ by the triangles $\langle w, v_{1}, v_{4} \rangle$, $\langle v_{1}, v_{i}, v_{i+1} \rangle, 2 \leq i \leq 3$.
We define $f'$ such that it coincides with $f$ on all simplices
which are common to $D$ and $D'$. We define $f'$ such that $f'(\langle v_{1}, v_{i} \rangle) = \langle f(v_{1}), f(v_{i}) \rangle$, $3 \leq i \leq 4$, $f'(\langle v_{1},v_{i},v_{i+1} \rangle) = \langle f(v_{1}),f(v_{i}),f(v_{i+1}) \rangle$, $2 \leq i \leq 3$, \\ $f'(\langle w,v_{1},v_{4} \rangle) = \langle f(w),f(v_{1}),f(v_{4}) \rangle$.
Because $f$ is simplicial, so is $f'$. Then $(D', f')$ is indeed a filling diagram for $\gamma$.
Note that $D$ and $D'$ have the same area. So $D'$ has minimal area. There is a full $4$-cycle $(v_{1},w,v_{7},v_{8})$ in $D'_{v}$, and a full $6$-cycle  $(v,v_{1},v_{4},v_{5},v_{6},v_{7})$ in $D'_{w}$. By case $3$, $k=4$ (or case $1$, $k=6$) we have reached a contradiction.

c. Suppose $f(v_{1}) \sim f(v_{5})$. According to Lemma \ref{3.2}, the $5$-cycle $(f(v_{1}),$ $f(v_{2}),f(v_{3}),f(v_{4}),f(v_{5}))$ in $X_{f(w)}$ is not full. So at least one of the following holds $f(v_{1}) \sim f(v_{i})$,  $3 \leq i \leq 4$. Assume w.l.o.g. $f(v_{1}) \sim f(v_{4})$. The other cases can be treated similarly. According to Lemma \ref{3.2}, the $4$-cycle $(f(v_{1}),f(v_{2}),f(v_{3}),f(v_{4}))$ in $X_{f(w)}$ is not full. Assume w.l.o.g. $f(v_{1}) \sim f(v_{3})$. We choose a filling diagram $(D', f')$ for $\gamma$ such that $D'$ is triangulated with the same simplices like $D$ except for the triangles $\langle w, v_{i}, v_{i+1} \rangle, 1 \leq i \leq 4$ in $D$ which are replaced in
$D'$ by the triangles $\langle w, v_{1}, v_{5} \rangle$, $\langle v_{1}, v_{i}, v_{i+1} \rangle, 2 \leq i \leq 4$.
We define $f'$ such that it coincides with $f$ on all simplices
which are common to $D$ and $D'$. We define $f'$ such that $f'(\langle v_{1}, v_{i} \rangle) = \langle f(v_{1}), f(v_{i}) \rangle $, $3 \leq i \leq 5$, $f'(\langle v_{1},v_{i},v_{i+1} \rangle) = \langle f(v_{1}),f(v_{i}),f(v_{i+1}) \rangle$, $2 \leq i \leq 4$, $f'(\langle w,v_{1},v_{5} \rangle) = \langle f(w),f(v_{1}),f(v_{5}) \rangle$.
Because $f$ is simplicial, so is $f'$. Then $(D', f')$ is indeed a filling diagram for $\gamma$.
Note that the discs $D$ and $D'$ have the same area. So $D'$ has minimal area. There is a full $4$-cycle $(v_{1},w,v_{7},v_{8})$ in $D'_{v}$, and a full $5$-cycle  $(v,v_{1},v_{5},v_{6},v_{7})$ in $D'_{w}$. Hence case $2$, $k=4$ (or case $1$, $k=5$) implies a contradiction.

\begin{figure}[h]
    \begin{center}
       \includegraphics[height=4cm]{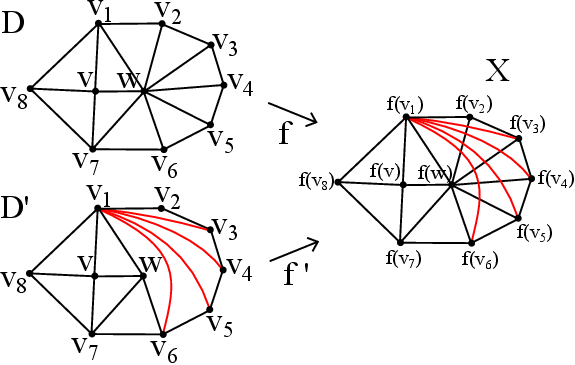}
      \caption{Case $5$, $f(v_{1}) \sim f(v_{6})$}
    \end{center}
\end{figure}

d. Suppose $f(v_{1}) \sim f(v_{6})$. According to Lemma \ref{3.2}, the $6$-cycle $(f(v_{1}),$ $f(v_{2}),f(v_{3}),f(v_{4}),f(v_{5}),f(v_{6}))$ in $X_{f(w)}$ is not full. So at least one of the following holds $f(v_{1}) \sim f(v_{i})$, $3 \leq i \leq 5$. Assume w.l.o.g. $f(v_{1}) \sim f(v_{5})$. The other cases can be treated similarly. According to Lemma \ref{3.2}, the $5$-cycle $(f(v_{1}),f(v_{2}),f(v_{3}),f(v_{4}),f(v_{5}))$ in $X_{f(w)}$ is not full. So at least one of the following holds $f(v_{1}) \sim f(v_{i})$, $3 \leq i \leq 4$. Assume w.l.o.g. $f(v_{1}) \sim f(v_{4})$. The other cases can be treated similarly. According to Lemma \ref{3.2}, the $4$-cycle $(f(v_{1}), f(v_{2}), f(v_{3}), f(v_{4}))$ in $X_{f(w)}$ is not full. Assume w.l.o.g. $f(v_{1}) \sim f(v_{3})$. We choose a filling diagram $(D', f')$ for $\gamma$ such that $D'$ is triangulated with the same simplices like $D$ except for the triangles $\langle w, v_{i}, v_{i+1} \rangle, 1 \leq i \leq 5$ in $D$ which are replaced in
$D'$ by the triangles $\langle w, v_{1}, v_{6} \rangle$, $\langle v_{1}, v_{i}, v_{i+1} \rangle, 2 \leq i \leq 5$.
We define $f'$ such that it coincides with $f$ on all simplices
which are common to $D$ and $D'$. We define $f'$ such that $f'(\langle v_{1}, v_{i} \rangle) = \langle f(v_{1}), f(v_{i}) \rangle$, $3 \leq i \leq 6$,  $f'(\langle v_{1},v_{i},v_{i+1} \rangle) = \langle f(v_{1}),f(v_{i}),f(v_{i+1}) \rangle$, $2 \leq i \leq 5$, $f'(\langle w,v_{1},v_{6} \rangle) = \langle f(w),f(v_{1}),f(v_{6}) \rangle$.
Because $f$ is simplicial, so is $f'$. Then $(D', f')$ is indeed a filling diagram for $\gamma$.
Because $D$ and $D'$ have the same area, $D'$ has minimal area. There is a full $4$-cycle $(v_{1},w,v_{7},v_{8})$ in $D'_{v}$, and a full $4$-cycle $(v,v_{1},v_{6},v_{7})$ in $D'_{w}$. Then case $1$, $k=4$ implies a contradiction.

e. Suppose $f(v_{1}) \sim f(v_{7})$. We consider a filling diagram $(D',f')$ for $\gamma$ such that $D'$ is triangulated with the same simplices like $D$ except for the triangles $\langle v,w,v_{1} \rangle$, $\langle v,w,v_{7} \rangle$, $\langle v,v_{7},v_{8} \rangle$, $\langle v,v_{1},v_{8} \rangle$ in $D$ which are replaced in $D'$ by the triangles $\langle w,v_{1},v_{7} \rangle$, $\langle v_{1},v_{7},v_{8} \rangle$. We define $f'$ such that it coincides with $f$ on all simplices
which are common to $D$ and $D'$. We define $f'$ such that $f'(\langle v_{1},v_{7} \rangle) = \langle f(v_{1}),f(v_{7}) \rangle$, $f'(\langle w,v_{1},v_{7} \rangle) = \langle f(w),f(v_{1}),f(v_{7}) \rangle$, $f'(\langle v_{1},v_{7},v_{8} \rangle) = \langle f(v_{1}),$ $f(v_{7}),f(v_{8}) \rangle$. Because $f$ is simplicial, $f'$ is also simplicial. So $(D',f')$ is a filling diagram for $\gamma$. Note that the area of $D'$ is less than the area of $D$. Hence we have reached a contradiction.

f. Suppose $f(v) \sim f(v_{2})$. We choose a filling diagram $(D', f')$ for $\gamma$ such that $D'$ is triangulated with the same simplices like $D$ except for the triangles $\langle w,v,v_{1} \rangle, \langle w,v_{1},v_{2} \rangle$  in $D$ which are replaced in
$D'$ by the triangles $\langle v, w, v_{2} \rangle$, $\langle v, v_{1}, v_{2} \rangle$.
We define $f'$ such that it coincides with $f$ on all simplices
which are common to $D$ and $D'$. We define $f'$ such that $f'(\langle v, v_{2} \rangle) = \langle f(v), f(v_{2}) \rangle$, $f'(\langle v,w,v_{2} \rangle) = \langle f(v),f(w),f(v_{2}) \rangle$, $f'(\langle v,v_{1},v_{2} \rangle) $ $= \langle f(v),$ $f(v_{1}),f(v_{2}) \rangle$.
Because $f$ is simplicial, so is $f'$. Then $(D', f')$ is indeed a filling diagram for $\gamma$.
Because $D$ and $D'$ have the same area, $D'$ has minimal area. There is a full $5$-cycle $(v_{1},v_{2},w,v_{7},v_{8})$ in $D'_{v}$, and a full $7$-cycle  $(v,v_{2},v_{3},v_{4},v_{5},v_{6},v_{7})$ in $D'_{w}$. Then case $4$, $k=5$ implies a contradiction.

g. Suppose $f(v) \sim f(v_{3})$. Note there is a $4$-cycle $\alpha = (f(v),f(v_{1}),$ $f(v_{2}),f(v_{3}))$ in $X_{f(w)}$ and there is a $4$-cycle $(f(v_{1}),f(w),$ $f(v_{7}),f(v_{8}))$ in $X_{f(v)}$. Lemma \ref{3.2} implies that $\alpha$ is not full. So either $f(v) \sim f(v_{2})$ or $f(v_{1}) \sim f(v_{3})$. Both cases are treated below.

$\bullet$ Suppose $f(v) \sim f(v_{2})$. We choose a filling diagram $(D', f')$ for $\gamma$ such that $D'$ is triangulated with the same simplices like $D$ except for the triangles $\langle w,v,v_{1} \rangle, \langle w,v_{i},v_{i+1} \rangle, 1 \leq i \leq 2$ in $D$ which are replaced in
$D'$ by the triangles $\langle v, w, v_{3} \rangle,$ $\langle v, v_{i}, v_{i+1} \rangle, 1 \leq i \leq 2$.
We define $f'$ such that it coincides with $f$ on all simplices
which are common to $D$ and $D'$. We define $f'$ such that $f'(\langle v, v_{i} \rangle) = \langle f(v), f(v_{i}) \rangle$, $2 \leq i \leq 3$, $f'(\langle v,v_{i},v_{i+1} \rangle) = \langle f(v),f(v_{i}),f(v_{i+1}) \rangle$, $1 \leq i \leq 2$, $f'(\langle v,w,v_{3} \rangle) = \langle f(v),f(w),f(v_{3}) \rangle$.
Because $f$ is simplicial, so is $f'$. Then $(D', f')$ is indeed a filling diagram for $\gamma$.
Because $D$ and $D'$ have the same area, $D'$ has minimal area. There is a full $6$-cycle $(v_{1},v_{2},v_{3},w,v_{7},v_{8})$ in $D'_{v}$, and a full $6$-cycle  $(v,v_{3},v_{4},v_{5},v_{6},v_{7})$ in $D'_{w}$. By case $3$, $k=6$ we have reached a contradiction.

$\bullet$ Suppose $f(v_{1}) \sim f(v_{3})$. We choose a filling diagram $(D', f')$ for $\gamma$ such that $D'$ is triangulated with the same simplices like $D$ except for the triangles $\langle w,v,v_{1} \rangle, \langle w,v_{i},v_{i+1} \rangle, 1 \leq i \leq 2$ in $D$ which are replaced in
$D'$ by the triangles $\langle v, w, v_{3} \rangle,$ $\langle v, v_{1}, v_{3} \rangle$, $\langle v_{1}, v_{2}, v_{3} \rangle$.
We define $f'$ such that it coincides with $f$ on all simplices
which are common to $D$ and $D'$. We define $f'$ such that $f'(\langle v, v_{3} \rangle) = \langle f(v), f(v_{3}) \rangle,$ $f'(\langle v_{1}, v_{3} \rangle) = \langle f(v_{1}), f(v_{3}) \rangle$, $f'(\langle v,w,v_{3} \rangle) = \langle f(v),f(w),f(v_{3}) \rangle$, $f'(\langle v,v_{1},v_{3} \rangle) = \langle f(v),f(v_{1}),f(v_{3}) \rangle$, $f'(\langle v_{1},v_{2},v_{3} \rangle) = \langle f(v_{1}),f(v_{2}),f(v_{3}) \rangle$.
Because $f$ is simplicial, so is $f'$. Then $(D', f')$ is indeed a filling diagram for $\gamma$.
Because $D$ and $D'$ have the same area, $D'$ has minimal area. There is a full $5$-cycle $(v_{1},v_{3},w,v_{7},v_{8})$ in $D'_{v}$, and a full $6$-cycle  $(v,v_{3},v_{4},v_{5},v_{6},v_{7})$ in $D'_{w}$. Hence case $2$, $k=6$ (or case $3$, $k=5$) implies a contradiction.

\begin{figure}[h]
    \begin{center}
       \includegraphics[height=4cm]{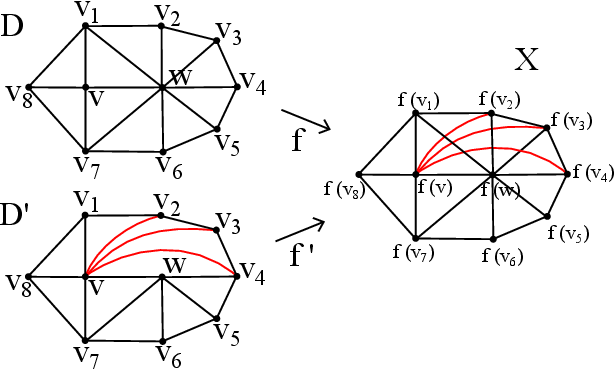}
      \caption{Case $5$, $f(v) \sim f(v_{4})$}
    \end{center}
\end{figure}

h. Suppose $f(v) \sim f(v_{4})$. Note that there is a $5$-cycle $\alpha = (f(v), f(v_{1}),$ $f(v_{2}),$ $ f(v_{3}), f(v_{4}))$ in $X_{f(w)}$ and there is a $4$-cycle $(f(v_{1}), f(w),$ $f(v_{7}),$ $ f(v_{8}))$ in $X_{f(v)}$. Lemma \ref{3.2} implies that $\alpha$ is not full. So at least one of the following holds $f(v) \sim f(v_{i})$, $2 \leq i \leq 3$. Assume w.l.o.g. $f(v) \sim f(v_{3})$. The other cases can be treated similarly. Lemma \ref{3.2} implies that the $4$-cycle $(f(v), f(v_{1}), f(v_{2}), f(v_{3}))$ in $X_{f(w)}$ is not full. Assume w.l.o.g. $f(v) \sim f(v_{2})$. We choose a filling diagram $(D', f')$ for $\gamma$ such that $D'$ is triangulated with the same simplices like $D$ except for the triangles $\langle w,v,v_{1} \rangle, \langle w,v_{i},v_{i+1} \rangle, 1 \leq i \leq 3$ in $D$ which are replaced in
$D'$ by the triangles $\langle v, w, v_{4} \rangle,$ $\langle v, v_{i}, v_{i+1} \rangle, 1 \leq i \leq 3$. We define $f'$ such that it coincides with $f$ on all simplices
which are common to $D$ and $D'$. We define $f'$ such that $f'(\langle v, v_{i} \rangle) = \langle f(v), f(v_{i}) \rangle$, $2 \leq i \leq 4$, $f'(\langle v,v_{i},v_{i+1} \rangle) = \langle f(v),f(v_{i}),f(v_{i+1}) \rangle$, $1 \leq i \leq 3$, $f'(\langle v,w,v_{4} \rangle) = \langle f(v), f(w), f(v_{4})\rangle$.
Because $f$ is simplicial, so is $f'$. Then $(D', f')$ is indeed a filling diagram for $\gamma$. Because $D$ and $D'$ have the same area, $D'$ has minimal area. There is a full $5$-cycle $(v,v_{4},v_{5},v_{6},v_{7})$ in $D'_{w}$, and a full $7$-cycle $(v_{1},v_{2},v_{3},v_{4},w,v_{7},v_{8})$ in $D'_{v}$. Then by case $4$, $k=5$ we have reached a contradiction.

Hence $w$ is $9$-large.

According to the $5$ cases above, $D$ satisfies the conditions $(5/9)$, $(6/8)$ and $(7/7)$ in the definition of a $5/9$-complex.
Thus $D$ is a $5/9$-complex.

\end{enumerate}
\end{proof}

\section{$5/9$-disc diagrams}

In this section we study properties of $5/9$-disc diagrams.

\begin{lemma}\label{4.1}
Let $X$ be a $5/9$-complex and let $\beta$ be a homotopically trivial loop in $X$. Let $(D,f)$ be a minimal filling diagram for $\beta$. Let $v$ and $w$ be adjacent interior vertices of $D$ such that $v$ is $4$-large but not $5$-large and $w$ is $9$-large. Let $D_{v} = (v_{1}, w, v_{8}, v_{9})$, let $D_{w} = (v_{1}, v_{2}, v_{3}, v_{4}, v_{5}, v_{6}, v_{7}, v_{8}, v)$ and let $\gamma = (v_{1}, v_{2}, v_{3}, v_{4}, v_{5}, v_{6}, v_{7}, v_{8}, v_{9})$. Then:

a. $\gamma$ is full in $D$;

b. $\gamma$ is not contained in the link of a vertex.

\end{lemma}

\begin{proof}

Because $X$ is a $5/9$-complex and $D$ has minimal area, Proposition \ref{3.3} implies that $D$ is a $5/9$-complex.

a. Suppose that $\gamma$ is not full in $D$. So $\gamma$ has at least one diagonal $\langle v_{1},v_{j} \rangle, 3 \leq j \leq 8$. Since $D$ is flat, such a diagonal lies outside both of $D_{w}$ and of $D_{v}$. Because the map $f$ is simplicial and nondegenerate, there is an edge $\langle f(v_{1}),$ $f(v_{j}) \rangle, 3 \leq j \leq 8$ in $X$. In the six cases below we will show, by reaching each time contradiction, that $\gamma$ has no diagonal $\langle v_{1},v_{j} \rangle$, $3 \leq j \leq 8$.

 \begin{enumerate}

 \item Case $1$: $j=3$.

We consider a filling diagram $(D',f')$ for $\beta$. The disc $D'$ is triangulated with the same simplices like $D$ except for the triangles $\langle w, v_{i}, v_{i+1} \rangle, 1 \leq i \leq 2$ in $D$ which are replaced in $D'$ by the triangles $\langle w, v_{1}, v_{3} \rangle, \langle v_{1},v_{2},v_{3} \rangle$. We define $f'$ such that it coincides with $f$ on all simplices which are common to $D$ and $D'$. We define $f'$ such that $f'(\langle v_{1},v_{3} \rangle) = \langle f(v_{1}), f(v_{3})\rangle, f'(\langle w,v_{1},v_{3} \rangle)$ $= \langle f(w),f(v_{1}),f(v_{3}) \rangle, f'(\langle v_{1},v_{2},v_{3} \rangle) = \langle f(v_{1}),$ $f(v_{2}),$ $f(v_{3}) \rangle$. Hence, because $f$ is simplicial, $f'$ is also simplicial. Thus $(D',f')$ is indeed a filling diagram for $\beta$. Because the discs $D$ and $D'$ have the same area, $D'$ has minimal area. Then Proposition \ref{3.3} implies that $D'$ is a $5/9$-complex. Note that in $D'$, $v$ is $4$-large but not $5$-large and $w$ is $8$-large. Since $v$ and $w$ are adjacent, this implies a contradiction with the $(5/9)$-condition satisfied by $D'$.

\item Case $2$: $j=4$.

Note that there is a $4$-cycle $(f(v_{1}), f(w), f(v_{8}), f(v_{9}))$ in $X_{f(v)}$ and there is a $4$-cycle $\alpha = (f(v_{1}), f(v_{2}), f(v_{3}), f(v_{4}))$ in $X_{f(w)}$. Lemma \ref{3.2} implies that $\alpha$ is not full. Assume w.l.o.g. $f(v_{1}) \sim f(v_{3})$.
We consider a filling diagram $(D',f')$ for $\beta$ such that $D'$ is triangulated with the same simplices like $D$ except for the triangles $\langle w, v_{i}, v_{i+1} \rangle, 1 \leq i \leq 3$ in $D$ which are replaced in $D'$ by the triangles $\langle w, v_{1}, v_{4} \rangle,$ $\langle v_{1}, v_{i},v_{i+1} \rangle, 2 \leq i \leq 3$. We define $f'$ such that it coincides with $f$ on all simplices which are common to $D$ and $D'$. We define $f'$ such that $f'(\langle v_{1},v_{i} \rangle) = \langle f(v_{1}), f(v_{i})\rangle,$ $3 \leq i \leq 4$, $f'(\langle v_{1},v_{i},v_{i+1} \rangle) =$ $\langle f(v_{1}),f(v_{i}),f(v_{i+1}) \rangle, 2 \leq i \leq 3$, $f'(\langle w,v_{1},v_{4} \rangle) = \langle f(w),f(v_{1}),f(v_{4}) \rangle$. Because $f$ is simplicial, $f'$ is also simplicial. Thus $(D',f')$ is indeed a filling diagram for $\beta$. Because $D$ and $D'$ have the same area, $D'$ has minimal area. Then Proposition \ref{3.3} implies that $D'$ is a $5/9$-complex. Note that in $D'$, $v$ is $4$-large but not $5$-large and $w$ is $7$-large. Hence we have reached a contradiction.

\item Case $3$: $j=5$.

Note that there is a $5$-cycle $(f(v_{1}), f(v_{2}), f(v_{3}), f(v_{4}), f(v_{5}))$ in $X_{f(w)}$. Lemma \ref{3.2} implies that this $5$-cycle is not full. So at least one of the following holds $f(v_{1}) \sim f(v_{i}), 3 \leq i \leq 4$. Assume w.l.o.g. $f(v_{1}) \sim f(v_{4})$. The other cases can be treated similarly. Note that there is a $4$-cycle $(f(v_{1}), f(v_{2}), f(v_{3}), f(v_{4}))$ in $X_{f(w)}$ which, according to Lemma \ref{3.2}, is not full. Assume w.l.o.g. $f(v_{1}) \sim f(v_{3})$.
We consider a filling diagram $(D',f')$ for $\beta$ such that $D'$ is triangulated with the same simplices like $D$ except for the triangles $\langle w, v_{i}, v_{i+1} \rangle, 1 \leq i \leq 4$ in $D$ which are replaced in $D'$ by the triangles $\langle w, v_{1}, v_{5} \rangle, \langle v_{1},v_{i},v_{i+1} \rangle, 2 \leq i \leq 4$. We define $f'$ such that it coincides with $f$ on all simplices which are common to $D$ and $D'$. We define $f'$ such that $f'(\langle v_{1},v_{i} \rangle) = \langle f(v_{1}), f(v_{i})\rangle, 3 \leq i \leq 5$, $f'(\langle v_{1},v_{i},v_{i+1} \rangle) = \langle f(v_{1}),f(v_{i}),f(v_{i+1}) \rangle, 2 \leq i \leq 4$, $f'(\langle w,v_{1},v_{5} \rangle) = \langle f(w),f(v_{1}),f(v_{5}) \rangle$. Because $f$ is simplicial, $f'$ is also simplicial. Thus $(D',f')$ is indeed a filling diagram for $\beta$. Because $D$ and $D'$ have the same area, $D'$ is minimal. Hence Proposition \ref{3.3} implies that $D'$ is a $5/9$-complex. Note that in $D'$, $v$ is $4$-large but not $5$-large, while $w$ is $6$-large. This implies a contradiction.

\item Case $4$: $j=6$.

Note that there is a $6$-cycle $(f(v_{1}), f(v_{2}), f(v_{3}), f(v_{4}), f(v_{5}), f(v_{6}))$ in $X_{f(w)}$.  By Lemma \ref{3.2}, this $6$-cycle is not full. So at least one of the following holds $f(v_{1}) \sim f(v_{i})$, $3 \leq i \leq 5$. Assume w.l.o.g. $f(v_{1}) \sim f(v_{5})$. The other cases can be treated similarly. Note that there is a $5$-cycle $(f(v_{1}), f(v_{2}), f(v_{3}), f(v_{4}), f(v_{5}))$ in $X_{f(w)}$ which, according to Lemma \ref{3.2}, is not full. So at least one of the following holds $f(v_{1}) \sim f(v_{i})$, $3 \leq i \leq 4$. Assume w.l.o.g. $f(v_{1}) \sim f(v_{4})$. The other cases can be treated similarly. Note that there is a $4$-cycle $(f(v_{1}), f(v_{2}), f(v_{3}), f(v_{4}))$ in $X_{f(w)}$ which, according to Lemma \ref{3.2}, is not full. Assume w.l.o.g. $f(v_{1}) \sim f(v_{3})$.
We consider a filling diagram $(D',f')$ for $\beta$ such that $D'$ is triangulated with the same simplices like $D$ except for the triangles $\langle w, v_{i}, v_{i+1} \rangle, 1 \leq i \leq 5$ in $D$ which are replaced in $D'$ by the triangles $\langle w, v_{1}, v_{6} \rangle, \langle v_{1},v_{i},v_{i+1} \rangle, 2 \leq i \leq 5$. We define $f'$ such that it coincides with $f$ on all simplices which are common to $D$ and $D'$. We define $f'$ such that $f'(\langle v_{1},v_{i} \rangle) = \langle f(v_{1}), f(v_{i})\rangle, 3 \leq i \leq 6$, $f'(\langle v_{1},v_{i},v_{i+1} \rangle) = \langle f(v_{1}),f(v_{i}),f(v_{i+1}) \rangle, 2 \leq i \leq 5$, $f'(\langle w,v_{1},v_{6} \rangle) = \langle f(w),f(v_{1}),f(v_{6}) \rangle$. Because $f$ is simplicial, $f'$ is also simplicial. Thus $(D',f')$ is indeed a filling diagram for $\beta$. Because $D$ and $D'$ have the same area, $D'$ is minimal. Proposition \ref{3.3} implies that $D'$ is a $5/9$-complex. Note that in $D'$, $v$ is $4$-large but not $5$-large, while $w$ is $5$-large. Hence we have reached a contradiction.

\item Case $5$: $j=7$.

Note that there is a $7$-cycle $(f(v_{1}), f(v_{2}), f(v_{3}), f(v_{4}), f(v_{5}), f(v_{6}), $ $f(v_{7}))$ in $X_{f(w)}$.  Lemma \ref{3.2} implies that this $7$-cycle is not full. So at least one of the following holds $f(v_{1}) \sim f(v_{i})$, $3 \leq i \leq 6$. Assume w.l.o.g. $f(v_{1}) \sim f(v_{6})$. The other cases can be treated similarly. Note that there is a $6$-cycle $(f(v_{1}), f(v_{2}), f(v_{3}), f(v_{4}), f(v_{5}), f(v_{6}))$ in $X_{f(w)}$ which, according to Lemma \ref{3.2}, is not full. So at least one of the following holds $f(v_{1}) \sim f(v_{i})$, $3 \leq i \leq 5$. Assume w.l.o.g. $f(v_{1}) \sim f(v_{5})$. The other cases can be treated similarly. Note that there is a $5$-cycle $(f(v_{1}), f(v_{2}), f(v_{3}), f(v_{4}), f(v_{5}))$ in $X_{f(w)}$ which, according to Lemma \ref{3.2}, is not full. So at least one of the following holds $f(v_{1}) \sim f(v_{i})$, $3 \leq i \leq 4$. Assume w.l.o.g. $f(v_{1}) \sim f(v_{4})$. Note that there is a $4$-cycle $(f(v_{1}), f(v_{2}), f(v_{3}), f(v_{4}))$ in $X_{f(w)}$ which, according to Lemma \ref{3.2}, is not full. Assume w.l.o.g. $f(v_{1}) \sim f(v_{3})$.
We consider a filling diagram $(D',f')$ for $\beta$ such that $D'$ is triangulated with the same simplices like $D$ except for the triangles $\langle w, v_{i}, v_{i+1} \rangle, 1 \leq i \leq 6$ in $D$ which are replaced in $D'$ by the triangles $\langle w, v_{1}, v_{7} \rangle, \langle v_{1},v_{i},v_{i+1} \rangle, 2 \leq i \leq 6$. We define $f'$ such that it coincides with $f$ on all simplices which are common to $D$ and $D'$. We define $f'$ such that $f'(\langle v_{1},v_{i} \rangle) = \langle f(v_{1}), f(v_{i})\rangle, 3 \leq i \leq 7$, $f'(\langle v_{1},v_{i},v_{i+1} \rangle) = \langle f(v_{1}),f(v_{i}),f(v_{i+1}) \rangle, 2 \leq i \leq 6$, $f'(\langle w,v_{1},v_{7} \rangle) = \langle f(w),f(v_{1}),f(v_{7}) \rangle$. Because $f$ is simplicial, $f'$ is also simplicial. Thus $(D',f')$ is indeed a filling diagram for $\beta$. Note that $D$ and $D'$ have the same area. Since $D'$ is minimal, Proposition \ref{3.3} implies that it is a $5/9$-complex. Note that in $D'$, $v$ is $4$-large but not $5$-large, while $w$ is $4$-large. This implies a contradiction.

\item Case $6$: $j=8$.

 We consider a filling diagram $(D',f')$ of $\gamma$ such that $D'$ is triangulated with the same simplices like $D$ except for the triangles $\langle v,w,v_{1} \rangle$, $\langle v,w,v_{8} \rangle$, $\langle v,v_{8},v_{9} \rangle$, $\langle v,v_{1},v_{9} \rangle$ in $D$ which are replaced in $D'$ by the triangles $\langle w,v_{1},v_{8} \rangle$, $\langle v_{1},v_{8},v_{9} \rangle$. We define $f'$ such that it coincides with $f$ on all simplices which are common to $D$ and $D'$. We define $f'$ such that $f'(\langle v_{1},v_{8} \rangle) = \langle f(v_{1}),f(v_{8}) \rangle$, $f'(\langle w,v_{1},v_{8} \rangle) = \langle f(w),f(v_{1}),f(v_{8}) \rangle$, \\ $f'(\langle v_{1},v_{8},v_{9} \rangle) = \langle f(v_{1}),f(v_{8}),f(v_{9}) \rangle$. Because $f$ is simplicial, $f'$ is also simplicial. So $(D',f')$ is a filling diagram for $\gamma$. Note that the area of $D'$ is less than the area of $D$. Because $D$ has minimal area, we have reached a contradiction.

The six cases above ensure that
$\gamma$ has no diagonal $\langle v_{1},v_{j} \rangle$, $3 \leq j \leq 8$. So $\gamma$ is full in $D$.

\end{enumerate}

b. Assume there is a vertex $z$ in $D$ such that $\gamma$ is contained in its link. Since $D$ is flat, we may w.l.o.g. consider the vertex $z$ such that the vertices $v_{i}, 1 \leq i \leq 9$ lie w.r.t. $z$ as in the figure below.

Suppose there are two edges $\langle z,v_{2} \rangle$ running outside the region of $D$ bounded by $\gamma$: one on the left side of $\gamma$ and the other on its right side. Then the intersection of these edges is two vertices, namely their common faces. But in a simplicial complex the intersection of two simplices is either a single common face of both or the empty set. Because we have reached a contradiction, there is only one edge $\langle z,v_{2} \rangle$ running outside of $\gamma$.

\begin{figure}[h]
    \begin{center}
       \includegraphics[height=3.5cm]{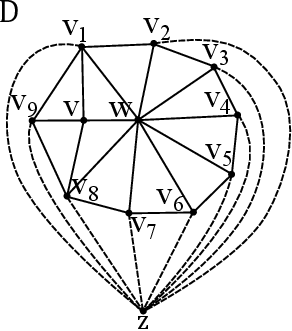}
      \caption{$\gamma$ is not contained in the link of a vertex}
    \end{center}
\end{figure}

Note that there is a full $4$-cycle $(v_{1},w,v_{8},v_{9})$ in $D_{v}$ and there is a full $4$-cycle $(v_{1},v,v_{8},z)$ in $D_{v_{9}}$. Hence the vertices $v$ and $v_{9}$ are both $4$-large. Because these vertices are adjacent, the $(5/9)$-condition fulfilled by $D$ implies a contradiction. So $\gamma$ is not contained in the link of a vertex.

\end{proof}

The proofs of the two lemmas below is similar to the proof of Lemma \ref{4.1}.

\begin{lemma}\label{4.2}
Let $X$ be a $5/9$-complex and let $\beta$ be a homotopically trivial loop in $X$. Let $(D,f)$ be a minimal filling diagram for $\beta$. Let $v$ and $w$ be adjacent interior vertices of $D$ such that $v$ is $5$-large but not $6$-large and $w$ is $8$-large. Let $D_{v} = (v_{1}, w, v_{7}, v_{8}, v_{9})$, let $D_{w} = (v_{1}, v_{2}, v_{3}, v_{4}, v_{5}, v_{6}, v_{7}, v)$ and let $\gamma = (v_{1}, v_{2}, v_{3}, v_{4}, v_{5}, v_{6}, v_{7}, v_{8}, v_{9})$. Then:

a. $\gamma$ is full in $D$;

b. $\gamma$ is not contained in the link of a vertex.

\end{lemma}

\begin{lemma}\label{4.3}
Let $X$ be a $5/9$-complex and let $\beta$ be a homotopically trivial loop in $X$. Let $(D,f)$ be a minimal filling diagram for $\beta$. Let $v$ and $w$ be adjacent interior vertices of $D$ such that $v$ is $6$-large but not $7$-large and $w$ is $7$-large. Let $D_{v} = (v_{1}, w, v_{6}, v_{7}, v_{8}, v_{9})$, let $D_{w} = (v_{1}, v_{2}, v_{3}, v_{4}, v_{5}, v_{6}, v)$ and let $\gamma = (v_{1}, v_{2}, v_{3}, v_{4}, v_{5}, v_{6}, v_{7}, v_{8}, v_{9})$. Then:

a. $\gamma$ is full in $D$;

b. $\gamma$ is not contained in the link of a vertex.

\end{lemma}

\section{The Gromov hyperbolicity of $5/9$-complexes}

In this section we give an application for the minimal filling diagrams lemma for $5/9$-complexes.

We start by showing that a minimal disc diagram associated to a homotopically trivial cycle in a $5/9$-complex, is $8$-located.

\begin{theorem}[$5/9$-diagrams are $8$-located]\label{5.1}
Let $X$ be a $5/9$-complex and let $\beta$ be a homotopically trivial loop in $X$. Let $(D,f)$ be a minimal filling diagram for $\beta$. Then $D$ is $8$-located.
\end{theorem}

\begin{proof}

Because $X$ is a $5/9$-complex and $D$ has minimal area, Proposition \ref{3.3} implies that $D$ is a $5/9$-complex.
Because $D$ is simply connected, any loop in $D$ is homotopically trivial.

Let $(l_{1},l_{2})$ be one of the pairs $\{(4,9),(5,8),(6,7)\}$.
Let $w_{1}$ be an interior vertex of $D$ which is $l_{1}$-large but not $(l_{1} + 1)$-large. Because $D$ is a $5/9$-complex, any neighbor of $w_{1}$
is $l_{2}$-large. Let $w_{2}$ be an interior vertex of $D$ adjacent to $w_{1}$. Let $D_{w_{1}} = (v_{1}, ..., v_{i}, w_{2}, v_{8}, v_{9})$, let $D_{w_{2}} = (v_{i}, ...,  v_{4}, v_{5}, v_{6}, v_{7}, v_{8}, w_{1})$, $1 \leq i \leq 3$ and let $\alpha =  (v_{1}, v_{2}, v_{3}, v_{4}, v_{5}, v_{6}, v_{7}, v_{8}, v_{9})$.
The Lemmas \ref{4.1}, \ref{4.2} and \ref{4.3} imply that $\alpha$ is full
and that it is not contained in the link of a vertex.
Because $\alpha$ is full, $|\alpha| = 9$.

\begin{figure}[h]
    \begin{center}
       \includegraphics[height=2cm]{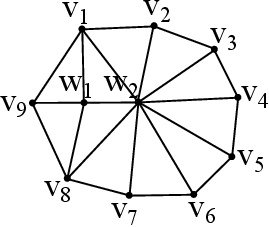}
      \caption{$|\alpha| = 9$, $w_{1}$ is $4$-large but not $5$-large, $w_{2}$ is $9$-large}
    \end{center}
\end{figure}

\begin{figure}[h]
    \begin{center}
       \includegraphics[height=2cm]{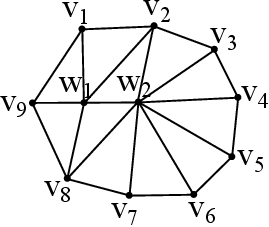}
      \caption{$|\alpha| = 9$, $w_{1}$ is $5$-large but not $6$-large, $w_{2}$ is $8$-large}
    \end{center}
\end{figure}

\begin{figure}[h]
    \begin{center}
       \includegraphics[height=2cm]{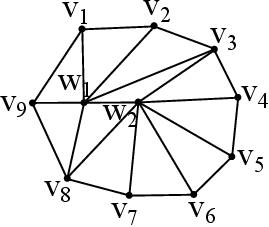}
      \caption{$|\alpha| = 9$, $w_{1}$ is $6$-large but not $7$-large, $w_{2}$ is $7$-large}
    \end{center}
\end{figure}

Suppose $D$ is not $8$-located. Then there exists a full loop $\gamma$ in $D$ of length at most $8$ which is not contained in the link of a vertex.
We consider at first the case when $|\gamma| = 8$.
Because $\gamma$ is full and homotopically trivial, there are at least two vertices in the interior of $\gamma$.
Assume at first there are two vertices $w_{1}, w_{2}$ inside $\gamma$ such that $w_{1}$ is $k_{1}$-large but not $(k_{1} + 1)$-large and $w_{2}$ is $k_{2}$-large. Let $\gamma_{j}$ be the cycle in $D_{w_{j}}, 1 \leq j \leq 2$. The flagness of $D$ implies that $k_{j} \geq 4, 1 \leq j \leq 2$. Because $D$ is flat and $|\gamma| = 8,$ we have $( k_{1},k_{2} ) \in \{ (4,8), (5,7), (6,6) \}$. Since $D$ is a $5/9$-complex, we have reached a contradiction.
We explain next why if $k_{1} = 4+i$, then $k_{2} = 8-i$, $0 \leq i \leq 2$. Because $w_{j}$ is an interior vertex of $\gamma$, there are two edges in the interior of $\gamma$ joining $w_{j}$ with two vertices of $\gamma$, $1 \leq j \leq 2$. So the $4+i$ edges of $\gamma_{1}$ are as follows: two of these edges lie in the interior of $\gamma$ and the other $2+i$ edges belong to $\gamma$, $0 \leq i \leq 2$. Since $|\gamma| = 8$, there are $6-i$ edges of $\gamma$ left which belong to $\gamma_{2}$, $0 \leq i \leq 2$. So the $8-i$ edges of $\gamma_{2}$ are as follows: two of these edges lie in the interior of $\gamma$ and the other $6-i$ edges belong to $\gamma$, $0 \leq i \leq 2$. Thus $k_{2} = 8-i$, $0 \leq i \leq 2$.

\begin{figure}[h]
    \begin{center}
       \includegraphics[height=1.7cm]{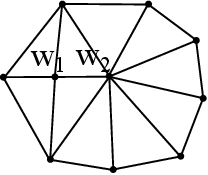}
         \label{}
      \caption{$|\gamma| = 8$, $w_{1}$ is $4$-large but not $5$-large, $w_{2}$ is $8$-large}
    \end{center}
\end{figure}

\begin{figure}[h]
    \begin{center}
       \includegraphics[height=1.7cm]{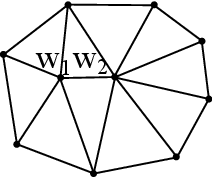}
        \label{}
      \caption{$|\gamma| = 8$, $w_{1}$ is $5$-large but not $6$-large, $w_{2}$ is $7$-large}
    \end{center}
\end{figure}


\begin{figure}[h]
    \begin{center}
       \includegraphics[height=1.7cm]{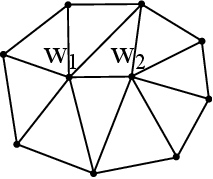}
        \label{}
      \caption{$|\gamma| = 8$, $w_{1}$ is $6$-large but not $7$-large, $w_{2}$ is $6$-large}
    \end{center}
\end{figure}


We consider further, besides the case $|\gamma| = 8$, the cases $4 \leq |\gamma| \leq 7$. Because $D$ is flat, only the following situations may occur.
If $|\gamma| = 7$, then $(k_{1},k_{2}) \in \{ (4,7), (5,6) \}$.
If $|\gamma| = 6$, then $(k_{1},k_{2}) \in \{ (4,6), (5,5) \}$.
If $|\gamma| = 5$, then $(k_{1},k_{2}) = (4,5)$.
If $|\gamma| = 4$, then $(k_{1},k_{2}) = (4,4)$.

 Besides the two vertices $w_{1}, w_{2}$ which already lie inside $\gamma$, we add inside $\gamma$ one more vertex $w_{3}$  which is $k_{3}$-large. Once we have added $w_{3}$ inside $\gamma$, the vertices $w_{i}$ become $k_{i}'$-large, $1 \leq i \leq 2$.
As explained below why, by adding $w_{3}$ inside $\gamma$, we return each time to a previous case reaching hereby a contradiction with the $5/9$-condition satisfied by $D$.


a. Assume $|\gamma| = 8$.

$\bullet$ Let $(k_{1},k_{2}) = (4,8)$. If we add $w_{3}$ inside $\gamma_{1}$, since $|\gamma_{1}| = 4$, then $(k_{1}',k_{3}) = (4,4)$. If we add $w_{3}$ inside $\gamma_{2}$ then, since $|\gamma_{2}| = 8$, $(k_{2}', k_{3}) \in \{ (4,8), (8,4), (5,7),$ $(7,5),$ $(6,6) \}$.

$\bullet$ Let $(k_{1},k_{2}) = (5,7)$. If we add $w_{3}$ inside $\gamma_{1}$, since $|\gamma_{1}| = 5$, then $(k_{1}',k_{3}) \in \{ (4,5),(5,4) \}$. If we add $w_{3}$ inside $\gamma_{2}$ then, since $|\gamma_{2}| = 7$, $(k_{2}', k_{3}) \in \{ (4,7), (7,4),$ $(5,6),$ $(6,5) \}$.

$\bullet$ Let $(k_{1},k_{2}) = (6,6)$. If we add $w_{3}$ inside $\gamma_{i}$, since $|\gamma_{i}| = 6$, then $(k_{i}', k_{3}) \in \{(4,6), (6,4),$ $(5,5)\}$, $1 \leq i \leq 2$.

b. Assume $|\gamma| = 7$.

$\bullet$ Let $(k_{1},k_{2}) = (4,7)$. If we add $w_{3}$ inside $\gamma_{1}$, since $|\gamma_{1}| = 4$, then $(k_{1}',k_{3}) = (4,4)$.  If we add $w_{3}$ inside $\gamma_{2}$ then, since $|\gamma_{2}| = 7$, $(k_{2}', k_{3}) \in \{ (4,7), (7,4),$ $(5,6),$ $(6,5) \}$.

$\bullet$ Let $(k_{1},k_{2}) = (5,6)$. If we add $w_{3}$ inside $\gamma_{1}$, since $|\gamma_{1}| = 5$, then $(k_{1}', k_{3}) \in \{(4,5), (5,4)\}$. If we add $w_{3}$ inside $\gamma_{2}$ then, since $|\gamma_{2}| = 6$, $(k_{2}', k_{3}) \in \{ (4,6),$ $(6,4),$ $(5,5) \}$.

c. Assume $|\gamma| = 6$.

$\bullet$ Let $(k_{1},k_{2}) \in \{(4,6)\}$. If we add $w_{3}$ inside $\gamma_{1}$, since $|\gamma_{1}| = 4$, then $(k_{1}', k_{3}) = (4,4)$. If we add $w_{3}$ inside $\gamma_{2}$ then, since $|\gamma_{2}| = 6$, $(k_{2}', k_{3}) \in \{ (4,6), (6,4), (5,5) \}$.

$\bullet$ Let $(k_{1},k_{2}) = (5,5)$. If we add $w_{3}$ inside $\gamma_{i}$, since $|\gamma_{i}| = 5$, then $(k_{i}', k_{3}) \in \{ (4,5), (5,4) \}$, $1 \leq i \leq 2$.

d. Assume $|\gamma| = 5$.

$\bullet$ Let $(k_{1},k_{2}) = (4,5)$. If we add $w_{3}$ inside $\gamma_{1}$, since $|\gamma_{1}| = 4$, then $(k_{1}', k_{3}) = (4,4)$. If we add $w_{3}$ inside $\gamma_{2}$ then, since $|\gamma_{2}| = 5$, $(k_{2}', k_{3}) \in \{ (4,5),(5,4) \}$.

e. Assume $|\gamma| = 4$.

$\bullet$ Let $(k_{1},k_{2}) = (4,4)$. If we add $w_{3}$ inside $\gamma_{i}$, since $|\gamma_{i}| = 4$, then $(k_{i}', k_{3}) = (4,4)$, $1 \leq i \leq 2$.

Assuming there are more than three vertices inside $\gamma$, we get similarly contradiction with the $5/9$-condition fulfilled by $D$ by returning each time to a previous case.
Any full loop in $D$ of length at most $8$ is therefore contained in the link of a vertex. Thus $D$ is $8$-located.

\end{proof}

We give further an application for the fact that simply connected, $8$-located simplicial complexes are Gromov hyperbolic.

\begin{corollary}[$5/9$-complexes are $8$-located]\label{5.2}
Let $X$ be a simply connected $5/9$-complex. Then $X$ is $8$-located. In particular, $X^{(0)}$ equipped with a path metric induced from $X^{(1)}$ is $\delta$-hyperbolic for a universal
constant $\delta$.
\end{corollary}

\begin{proof}

Let $\alpha$ be a cycle in $X$. Let $(D, f)$ be a minimal filling diagram for $\alpha$. Because $D$ is simply connected, any loop in $D$ is homotopically trivial.
According to Proposition \ref{3.3}, $D$ is a $5/9$-complex. Then Theorem \ref{5.1}
 implies that $D$ is $8$-located. So any full loop $\gamma$ in $D$ of length at most $8$, is contained in the link of a vertex $v$ of $D$.
Since $f$ is simplicial and nondegenerate, there exists a loop $\gamma'$ in $X$ such that for each edge $e \subset \gamma$, there is an edge $e' \subset \gamma'$ such that $f(e) = e'$. So $|\gamma| = |\gamma'|$.
Also there is a vertex $v'$ of $X$ adjacent to all vertices of $\gamma'$ such that $f(v) = v'$. Lemma \ref{3.1} implies that the vertex $f(v)$ differs from any vertex in $X_{f(v)}$ and the map $f$ is injective on $X_{f(v)}$.
So any full loop in $X$ of length at most $8$ is contained in the link of a vertex.
 Therefore, since $X$ is flag,
it is $8$-located. Because $X$ is simply connected and $8$-located, \cite{L-8loc}, Theorem $4.3$ implies that it is
Gromov hyperbolic.
\end{proof}

\end{document}